\documentclass[12pt]{article}
  \usepackage[utf8]{inputenc}
 \usepackage{graphicx}
 \usepackage{latexsym}
 \usepackage{amsmath}
 \usepackage{verbatim}
 \usepackage{amssymb}
 \usepackage{mathrsfs}
 \usepackage{amsmath,amssymb}
 \DeclareMathAlphabet      {\mathbf}{OT1}{cmr}{bx}{n}
 \DeclareFontFamily{OT1}{pzc}{}
 \DeclareFontShape{OT1}{pzc}{m}{it}%
 {<-> s * [1.15] pzcmi7t}{}
 \DeclareMathAlphabet{\mathpzc}{OT1}{pzc}{m}{it}
 \usepackage{enumitem}
 \usepackage{ dsfont }
 \usepackage{ mathrsfs }

 \usepackage{amsbsy }
 \usepackage{graphics}

 \usepackage{amsthm}
 \newtheorem{thm}{Theorem}
 \newtheorem{theorem}{Theorem}[section]

 \newtheorem{lemma}[theorem]{Lemma}

 \newtheorem{definition}[theorem]{Definition}

 \newtheorem{prop}[theorem]{Proposition}

 \newtheorem{remark}{Remark}[section]

 \usepackage{romannum}
 \usepackage[toc,page,title,titletoc,header]{appendix}
 \usepackage{hyperref}
 \usepackage{amscd}
 \usepackage[all,cmtip]{xy}

 \usepackage{abstract}

\begin{document}   \small
 	\pagenumbering{arabic}

	\pagenumbering{arabic}
	\title{Symplectic embeddings of balls into disk prequantization bundles}
	\author{{ Guanheng Chen}}
	\date{}
	%\begin{titlingpage}
	\maketitle
	%Department of Mathematics, University of Adelaide \\
	\thispagestyle{empty}
	\begin{abstract}
A prequantization bundle is a circle bundle over a symplectic surface with negative Euler class. A connection 1-form induces a natural contact form on such a bundle. The purpose of this paper is to study symplectic embeddings of balls into the associated disk bundles. To this end, we compute the ECH cobordism maps induced by these disk bundles. In addition, we compute the ECH spectrum of the prequantization bundles over the sphere and the torus.

	\end{abstract}

\section{Introduction and Main results}

Let $(X_1, \omega_1)$ and $(X_2, \omega_2)$ be two compact symplectic manifolds of the same dimension. A symplectic embedding of $(X_1, \omega_1)$ into $(X_2, \omega_2)$ is a smooth embedding $\varphi$ such that $\varphi^*\omega_2 = \omega_1$. When $(X_1, \omega_1)$ and $(X_2, \omega_2)$ are toric domains in $\mathbb{R}^4$, there are many results have been achieved  by M. Hutchings, D. Cristofaro-Gardiner, K. Choi, D. Frenkel, and V. G. B.  Ramos \cite{H3, H5,  DCG, CCFHV}.
 Beyond the toric cases, B. Ferreira and V. G. B. Ramos \cite{FRV}, along with A. Vicente \cite{FRV2}, recently studied the embedding problems from balls to the
unit disk cotangent bundle of the sphere or real projective space.

Let $Y$ be a closed three-manifold equipped with a contact form $\lambda$ such that $\lambda \wedge d\lambda > 0$.  To study the symplectic embedding phenomena in four dimension, M. Hutchings introduces a sequence of numerical invariants:
\begin{equation*}
0 < c_1(Y,\lambda) \le c_2(Y,\lambda) \le c_3(Y,\lambda) \le ... \le  \infty.
\end{equation*}
associated to $(Y,\lambda)$ which he calls the \textbf{ECH spectrum} \cite{H3}. The results mentioned in the first paragraph more or less rely on the computations of the ECH spectrum.

The purpose of this paper is to study symplectic embeddings a disjoint union of Darboux balls into a prequantization disk bundle.  Furthermore, we give some computations of ECH spectrum for the prequantization bundles over the sphere or torus.

Roughly speaking, a prequantization bundle is a circle bundle over a symplectic surface with negative Euler class. A holomorphic curve in its symplecticization has certain $S^1$-symmetry due to the fibration structure. Based on  this  observation,  J. Nelson and M. Weiler compute the embedded contact homology of the prequantization bundles \cite{NW} (based on D. Farris's PhD thesis \cite{Fa}). Their computations play a crucial  role in our proof.

The precise definition of the prequantization bundles is as follows.  Let $(\Sigma, \omega_\Sigma)$ be a closed surface with a volume form such that $[\omega_\Sigma] \in H^2(\Sigma,\mathbb{R}) \cap H^2(\Sigma,\mathbb{Z})$ is integral. Let $\pi_E: E \to \Sigma$ be a complex line bundle with $c_1(E) = -[\omega_\Sigma]$.  Then $E$ is called a \textbf{prequantization line bundle}.
Let $e :=\langle  c_1(E), [\Sigma] \rangle $ denote the degree of $E$.

Fix a Hermitian metric $h$ and Hermitian connection 1-form $A_\nabla$ such that $\frac{i}{2\pi}F_{A_\nabla} = -\omega_\Sigma$,  where $F_{A_{\nabla}}$ is the curvature of $A_{\nabla}$.  This gives us a global angular form  $\alpha_\nabla \in \Omega^1(E\setminus \Sigma, \mathbb{R})$.   Under a unitary trivialization $U \times \mathbb{C}$, $\alpha_{\nabla}$ is of the form  $\frac{1}{2\pi}(d\theta -iA_{\nabla} \vert_{U})$, where $d\theta$ is the angular form of $\mathbb{C}$ and  $A_{\nabla} \vert_{U}$ is a $i\mathbb{R}$ valued 1-form. Therefore, we have $d\alpha_{\nabla} = \pi_E^*\omega_{\Sigma}$ over $E \setminus \Sigma$.
A natural symplectic form on $E$ is defined by
\begin{equation*}
\Omega := \pi_E^*\omega_\Sigma + d(\rho^2 \alpha_\nabla),
\end{equation*}
where $\rho$ is the radius coordinate of $E$ defined by the metric $h$. Extend $\Omega$ over the zero section $\Sigma$ by
\begin{equation*}
  d(\rho^2 \alpha_{\nabla}) \vert_{fiber}:= (\mbox{area form of } \mathbb{C}) / \pi \mbox{ and } d(\rho^2 \alpha_\nabla)(T\Sigma, \cdot):=0.
\end{equation*}

Let  $\pi: Y: = \{ \rho = 1 \} \to \Sigma$ be the unit circle subbundle of $E$. Since
\begin{equation*}
  \Omega=2\rho d\rho\wedge \alpha_{\nabla} + (\rho^2  +1)d\alpha_{\nabla}
\end{equation*}
away from $\Sigma$,  the Liouville vector field is $Z= \frac{1+ \rho^2}{2\rho^2} \rho \partial_{\rho}$. Hence, $\Omega $ induces a contact form
$ \lambda  =\Omega(Z, \cdot) =2\alpha_{\nabla}$  on  $Y$.  The  contact  manifold  $(Y, \lambda)$  is  called  the  \textbf{prequantization bundle} of $(\Sigma, \omega_{\Sigma})$. Let $DE := \{\rho \le 1\} $ be the unit disk subbundle of $E$. Then $(DE, \Omega)$
forms a natural symplectic filling of $(Y, \lambda)$.

Let $B^4(r) := \{ (z_1,z_2) \in \mathbb{C}^2 \mid \pi|z_1|^2 + \pi|z_2|^2 \leq r \}$ be a ball of radius $\sqrt{r/\pi}$ in $\mathbb{C}^2$.   It
 is equipped with the standard symplectic form $\omega_{\mathrm{std}}$ of $\mathbb{C}^2$. The \textbf{Gromov width} of a
 symplectic four-manifold $(X,\Omega_X)$ is defined by
 \begin{equation*}
c_{\mathrm{Gr}}(X,\Omega_X) := \sup\{ a \in \mathbb{R}_{\geq 0} \mid \exists\ \mbox{embedding }, \varphi: B^4(a) \to X, \ \varphi^*\Omega_X = \omega_{\mathrm{std}} \}.
 \end{equation*}

Our main results are as follows.
\begin{thm} \label{thm1}
Let $\varphi: \bigsqcup_{i=1}^k (B^4(r_i), \omega_{\mathrm{std}}) \to (DE,\Omega)$ be a symplectic embedding. Then $\sum_{i=1}^k r_i \leq k$. In particular, $c_{\mathrm{Gr}}(DE,\Omega) \leq 1$.
\end{thm}

We remark that for $k^2 < 2|e|$, the estimates in Theorem \ref{thm1} are better than the
 volume constraint.   Furthermore, they are also better than the constrains from the ECH spectrum (at least in the case that $\Sigma$ is the sphere or the torus).

The proof of Theorem \ref{thm1} relies on the computations of the ECH cobordism maps induced by $(DE, \Omega)$. From the construction, every fiber of  $\pi :Y \to \Sigma $ is a Reeb  orbit.   Since the Reeb orbits with fixed length are not isolated, the contact form $\lambda$  is
degenerate. To define the ECH group and the cobordism map, we follow \cite{NW} to perform a perturbation on $\lambda$ by a perfect Morse function $H:\Sigma \to \mathbb{R}$. Under this perturbation, the simple Reeb orbits with length less than $L_{\varepsilon}$ correspond one-to-one to the critical
points. Some suitable modifications also be made on $\Omega$. The results are denoted by $\lambda_{\varepsilon}$
and $\Omega_{\varepsilon}$ respectively. The details are given in Section \ref{section31}.

\paragraph{Coefficient} We use $\mathbb{F} :=\mathbb{Z}/2 \mathbb{Z}$-coefficient  throughout this paper.
\begin{thm} \label{thm2}
For any $0 < \varepsilon \ll 1$, let $\lambda_{\varepsilon}$ and $\Omega_\varepsilon$ be the perturbations of $\lambda$ and $\Omega$ defined in Section \ref{section31}, such that $\lambda_{\varepsilon}$ is $L_{\varepsilon}$-nondegenerate.  Fix $\Gamma \in \mathbb{Z}_{|e|}$. For any positive integer $M$ such that $\Gamma = M \mod |e|$ and $2M < L_\epsilon$, there exists $A \in H_2(DE, Y, \mathbb{Z})$ such that $\partial A =\Gamma$ and   the ECH cobordism map
\begin{equation*}
 {ECH}^{L_\varepsilon}(DE,\Omega_\varepsilon,A) : ECH^{L_{\varepsilon}}(Y, \lambda_{\varepsilon}, \Gamma) \to \mathbb{F}
\end{equation*}
 satisfies the following properties:
\begin{enumerate}
    \item $ECH^{L_\varepsilon}(DE,\Omega_\varepsilon,A)(e_{-}^M) = 1$, where $e_-$ is the Reeb orbit corresponding to the
 minimum of $H$;
    \item
Let  $\beta = e_+^{n_+} h_1^{n_1} \cdots h_{2g}^{n_{2g}} e_-^{n_-}$  be an ECH generator, where $h_i$ are the Reeb orbits
 corresponding to the saddle points of $H$, and $e_+$ is the Reeb orbit corresponding
 to the maximum of $H$.  Let  $N = n_+ + n_1 + \cdots + n_{2g} + n_-$.  If $N = M - d|e|$ for $d \geq 1$, then $ {ECH}^{L_\varepsilon}(DE,\Omega_\varepsilon,A)([\beta]) = 0$.
\end{enumerate}
\end{thm}
 Using Nelson and Weiler's computations on ECH, spectrality property of ECH
 spectrum, and the finiteness of ECH spectrum, we obtain the precise values of ECH spectrum in the following special cases.

\begin{thm} \label{thm3}
Suppose that $\Sigma$ is the two-sphere. Then for any $k \geq 0$, the $k$-th value on the ECH spectrum of $(Y,\lambda)$ is
\begin{equation*}
c_k(Y,\lambda) = 2d|e|,
\end{equation*}
where $d$ is the unique nonnegative integer such that
\begin{equation*}
2d + d|e|(d - 1) \leq 2k \leq 2d + d|e|(d + 1).
\end{equation*}
\end{thm}

\begin{thm} \label{thm4}
Suppose that  $\Sigma$ is the two-torus. Then we have the following statements.
\begin{itemize}
    \item
    If $k \neq \frac{n(n+1)|e|}{2}$ for any $n \in \mathbb{Z}_{\geq 1}$, then
 \begin{equation*}
   c_k(Y,\lambda) = 2 \lfloor \sqrt{\frac{2k}{|e|} + \frac{1}{4}} + \frac{1}{2} \rfloor |e|,
 \end{equation*}
 where $ \lfloor x \rfloor$ denote the maximal integer that is less than or equal to $x$;

    \item $c_{|e|}(Y,\lambda) = 4|e|$.
\end{itemize}
\end{thm}

\begin{remark}
If we consider $Y = \{\rho = c\}$, then the induced contact form on $Y$ is $\lambda_c = (1 + c^2)\alpha_{\nabla} = \frac{1 + c^2}{2} \lambda$.   By the conformality property of the ECH spectrum, i.e., $c_k(Y, \kappa \lambda) = \kappa c_k(Y, \lambda)$ for any positive constant $\kappa$ (see Remark 3.2 of \cite{H3}), we have
$$
c_k(Y, \lambda_c) = \frac{1 + c^2}{2} c_k(Y, \lambda).$$
\end{remark}

\begin{remark}
When $\Sigma$ is the sphere, there is a standard method to compute the $U$ map on $ECH^{L_\varepsilon} (Y, \lambda_\varepsilon)$ which  can be used to compute $c_k(Y, \lambda_\varepsilon)$. The methods are given by Hutchings when he computes the ECH of the three-sphere (Proposition 4.1 of \cite{H4}). In the results of Ferreira and Ramos \cite{FRV}, they use these arguments to compute the ECH spectrum of the cotangent circle bundle of the sphere.

We briefly sketch  the argument  in the case of prequantization bundles. By Nelson and Weiler's index computations (Proposition 3.5 of \cite{NW}), one can show that the holomorphic curves $\mathcal{C} $ contributing to the $U$ map are holomorphic cylinders. The degree of these cylinders is either zero or one. If the degree of $\mathcal{C} $ is zero, then it corresponds to an index 2 Morse flow line on $\mathbb{S}^2$. If the degree of $\mathcal{C} $ is one, then under the natural holomorphic structure of $E$, $\mathcal{C} $ is a
meromorphic section of $E$ with poles and zeros. Combining these facts, one could show that
\begin{enumerate}
  \item
  $U(e_+^i e_-^j ) = (e_+^{i-1}e_-^{j+1})$ (counting degree zero holomorphic cylinders);
  \item
$U(e_-^j ) = e_+^{j-|e|}$ (counting degree one holomorphic cylinders).
\end{enumerate}
 Here $e_{\pm}$ are defined in Section \ref{section31}.

\end{remark}

\begin{remark}
The computations of ECH cobordism maps can be used as an auxiliary tool to study the $U$ map on $ECH(Y, \lambda)$. Here is an example: Suppose that $e=-1$ and $\Sigma$ is the torus.  From the index computations of Nelson and Weiler and Lemma~\ref{lem12}, we obtain
\begin{equation*}
 U\begin{pmatrix}
e_-e_+ \\
h_1h_2
\end{pmatrix} = \begin{pmatrix}
a & 1 \\
b & c
\end{pmatrix}\begin{pmatrix}
e_+\\
e_-^2
\end{pmatrix},
\end{equation*}
 where  $a,b,c \in \mathbb{F}$ are some unknown numbers. By Theorem \ref{thm2}, we have  $c = ECH^L(DE, \Omega_\varepsilon, A)(U[h_1h_2])$. The methods in this paper can be used to show that the J-holomorphic currents $\mathcal{C}$ contributing to $ECH^L(DE, \Omega_\varepsilon, A)(U[h_1h_2])$ are embedded $J$-holomorphic curves passing
through a marked point. By Lemma \ref{lem13}, there exists  a $J$-holomorphic plane $\mathcal{C}'$ passing through the same marked point. However, due to  the homological reasons ($[\mathcal{C}] = Z_{h_1h_2} $ and $[\mathcal{C}']= Z_{e_+}$), these two curves should be disjoint. As a result, such a $\mathcal{C}$  cannot exist and $c = 0$. The advantage of the above argument is to avoid using the domain dependent almost complex structures as in \cite{NW} and all the methods are classical.

For computing the $U$ map, Farris \cite{Fa}, Nelson and Weiler (Page 11 of \cite{NW})  propose a heuristic idea about this. They suggest that the curves contributing to the $U$ map should be equivalent to the meromorphic multisections of $E$. We believe that to compute $ECH^L(DE, \Omega_\varepsilon, A) \circ U^k$  for a general class $A$ is also the same as counting the
meromorphic multisections of $E$. Because the cobordism maps are defined by counting
curves in $E$, their perspective appears more natural  in the cobordism setting.
\end{remark}

Now we sketch the idea of the proofs of the main results.
\paragraph{Idea of the proof}
The methods   of computing $ECH^L(DE, \Omega_{\varepsilon}, A )([\alpha_k])$ are more or less the same as \cite{GHC1}. Choose an  almost complex structure $J$ such that the fiber of $E$  over  the minimum of $H$ is $J$-holomorphic. For some fixed $A \in H_2(DE, Y, \mathbb{Z})$, we  show that a  cover  of the fiber  is  the only holomorphic current with $I=0$ in the
relative class  $A$. The result then follows from the  correspondence between solutions to the Seiberg-Witten equations and holomorphic curves (\cite{CG2}). The same
 argument can also be used to compute $ECH^{L_\varepsilon}(DE,\Omega_\varepsilon,A)(U^k[e^k_+])$. The non-vanishing
 of   $ECH^{L_\varepsilon}(DE,\Omega_\varepsilon,A)(U^k[e^k_+])$ yields   $J$-holomorphic curves passing through the
 centers of  $\varphi(\sqcup_{i=1}^k(B^4(r_i), \omega_{\mathrm{std}}))$. Then the standard monotonicity lemma gives the result
 of Theorem \ref{thm1}.

 When the base manifold is the sphere or the torus, thanks to the computations from P. Ozsv\'{a}th, Z. Szab\'{o} \cite{OS}, and K. Park \cite{Park},   the existence of $U$ sequences is known. Combining these results with the ``volume property" of the ECH spectrum \cite{CHR}, we
 obtain the finiteness of $c_k(Y,\lambda)$. Using spectrality, we derive   Theorem \ref{thm3} and
 the first item of Theorem \ref{thm4}. To obtain the second item of Theorem \ref{thm4}, we also need to compute $\langle U^ke^k_+,e^k_- \rangle$.

\section*{Acknowledgements}
The author sincerely thanks Nelson and Weiler, whose comments, suggestions, and corrections have greatly improved this paper. The author is also grateful for their selfless sharing of thoughts and ideas. Moreover, many parts of this paper rely significantly on their work. In addition, the author would like to thank the anonymous referee for his/her  helpful comments.

The author is supported by the National Natural Science Foundation of China
Youth Fund Project under Grant number: 12201106, and is partially supported by the National Natural Science Foundation of China under Grant number: 12271349.

\section{Preliminaries}
 In this section, we give a quick review of the embedded contact homology (abbreviated
 as ECH). For more details, please refer to  \cite{H4}.

Let $ Y$ be a closed, contact 3-dimensional manifold equipped with a nondegenerate contact
form $\lambda$. The contact structure of $Y$ is denoted by $\xi: =\ker \lambda$.
The \textbf{Reeb vector field} $R$ of $(Y, \lambda)$ is characterized by  conditions $\lambda(R)=1$ and $d\lambda(R, \cdot)=0$.
A \textbf{Reeb orbit} is a smooth map $\gamma: \mathbb{R}_{\tau} / T \mathbb{Z} \to Y $ satisfying the ODE $\partial_{\tau} \gamma =R \circ \gamma$ for some $T>0$. %The
%number $T$  is called the \textbf{action} of $\gamma$.

 The linearized Reeb flow for time $T$ defines a symplectic linear map
 \begin{equation*}
   P_{\gamma}: (\xi_{\gamma(0)}, d\lambda) \to (\xi_{\gamma(T)}, d\lambda).
 \end{equation*}
 The Reeb orbit $\gamma$ is called \textbf{nondegenerate} if $1$ is not an eigenvalue of $P_{\gamma}$.  A nondegenerate Reeb orbit $\gamma$ is called \textbf{elliptic} if the eigenvalues of $P_{\gamma}$ are on the unit circle, \textbf{positive hyperbolic} if the eigenvalues   are real positive
numbers, and \textbf{negative hyperbolic} if the eigenvalues are real negative numbers. Given $L \in \mathbb{R}$, the contact form $\lambda$ is called \textbf{$L$-nondegenerate} if all Reeb orbits with action less than $L$  are nondegenerate.  We  call $\lambda$  nondegenerate if all   Reeb orbits  are nondegenerate.

 An \textbf{orbit set}    $\alpha=\{(\alpha_i, m_i)\}$ is a finite set of Reeb orbits together with multiplicities, where $\alpha_i$ are distinct, nondegenerate, irreducible embedded Reeb orbits and $m_i$  are positive integers.  An orbit set $\alpha$ is called an \textbf{ECH generator} if $m_i=1$ whenever $\alpha_i$ is a hyperbolic orbit. Define the \textbf{action} of an orbit set $\alpha$ by
 \begin{equation*}
   \mathcal{A}_{\lambda}(\alpha):=\int_{\alpha} \lambda = \sum_{i} m_i \int_{\alpha_i} \lambda.
 \end{equation*}
 Note that  the action $   \int_{\gamma} \lambda$ is just the period $T$ in the definition of  Reeb orbits.

\textbf{Notation. }In the rest of the paper, we often write an orbit set using multiplicative notation $\alpha=\Pi_i \alpha_i^{m_i}$ instead.

% The intersection number $[\alpha] \cdot [fiber]$ is called the degree of $\alpha$.
The following definition will be used in the computations of the cobordism maps.
In fact, the elliptic Reeb orbits considered in our cases are one of the following two types of Reeb orbits.
\begin{definition} (see \cite{H5} Definition 4.1)
Fix  $L>0$.  Let  $\gamma$ be an embedded elliptic orbit with action $\mathcal{A}_{\lambda}(\gamma) < L$.
\begin{itemize}
\item
$\gamma$ is called $L$-positive elliptic if the rotation number $\theta  $ is in $ (0, \frac{\mathcal{A}_{\lambda}(\gamma) }{L}) \mod 1$.
\item
$\gamma$ is called $L$-negative elliptic if the rotation number $\theta $ is in $ ( -\frac{\mathcal{A}_{\lambda}(\gamma) }{L},0) \mod 1$.
\end{itemize}
\end{definition}

\paragraph{The ECH index}
Fix $\Gamma \in H_1(Y, \mathbb{Z})$. Given orbit sets $\alpha=\Pi_i  \alpha_i^{m_i}$ and $\beta=\Pi_j  \beta_j^{n_j}$ on  $Y$ with $[\alpha]=[\beta] =\Gamma$,  let $H_2(Y, \alpha ,\beta)$ be the set of 2-chains $Z$ such that $\partial Z = \sum_i m_i \alpha_i - \sum_j n_j \beta_j$, modulo boundaries of 3-chains.   An element in  $H_2(Y, \alpha ,\beta)$  is called \textbf{a relative homology class}. Note that the set  $H_2(Y, \alpha ,\beta)$ is an affine space over $H_2(Y, \mathbb{Z})$.

Given $Z \in H_2(Y, \alpha ,\beta)$ and  trivializations    $\tau$ of $\xi\vert_{\alpha}$ and $\xi \vert_{\beta}$,  the ECH index is defined by
\begin{equation*}
I(\alpha, \beta, Z) := c_{\tau}(\xi \vert_Z) + Q_{\tau}(Z) + \sum_i  \sum\limits_{p=1}^{m_i} CZ_{\tau}(\alpha_i^{p})- \sum_j \sum\limits_{q=1}^{n_j} CZ_{\tau}(\beta_j^{q}),
\end{equation*}
where $ c_{\tau}(\xi \vert_Z) $ and $Q_{\tau}(Z)$ are respectively the relative Chern number and the relative self-intersection number  (see \cite{H4} and \cite{H2}), and $CZ_{\tau}$ is the Conley-Zehnder index.  The  ECH index $I$  only depends  on orbit sets $\alpha$, $\beta$ and the  relative homology class  $Z$.

\paragraph{$J$-Holomorphic currents}
An almost complex structure on $(\mathbb{R} \times Y, d(e^s \lambda))$ is called \textbf{admissible} if $J$ is $\mathbb{R}$-invariant,  $J(\partial_s) =R$, $J(\xi) =\xi$ and $J \vert_{\xi}$ is  $d \lambda$-compatible. We denote set of admissible almost complex structures by $\mathcal{J}(Y, \lambda)$.

 A \textbf{$J$-holomorphic current} from $\alpha$ to $\beta$ is a formal sum $\mathcal{C} =\sum_a d_a C_a$, where $C_a$ are distinct  irreducible simple holomorphic curves, the $d_a$ are positive integers,  $\mathcal{C}$ is asymptotic to $\alpha$ as a current as $s \to \infty $ and  is asymptotic to $\beta$ as a current as $s \to -\infty.$ Fix $Z \in H_2(Y, \alpha, \beta)$. Let $\mathcal{M}^J(\alpha, \beta, Z)$ denote the moduli space of holomorphic currents with relative homology class $Z$.

  Let $C$ be $J$-holomorphic curve in $\mathbb{R} \times Y $ whose positive ends are asymptotic to $\alpha= \Pi \alpha_i^{m_i}$  and negative ends are asymptotic to $\beta= \Pi_j \beta_j^{n_j}$. For each $i$, let $k_i$ denotes the number of ends of $C$ at $\alpha_i$, and let $\{p_{ia}\}^{k_i}_{a=1}$ denote their multiplicities. Likewise,
for each $j$, let $l_j$ denote the number of ends of $u$ at $\beta_j$, and let $\{q_{jb}\}^{lj}_{b=1}$ denote their multiplicities.  The set of numbers $\{p_{ia}\}^{k_i}_{a=1}$ modulo order is  called the \textbf{partition} of $C$ at $\alpha_i$.
The \textbf{Fredholm index} of $u$ is defined by
\begin{eqnarray*}
\operatorname{ind} (C ):= -\chi(C) + 2 c_{\tau}(\xi \vert_{C} ) + \sum\limits_i \sum\limits_{a=1}^{k_i}CZ_{\tau} (\alpha^{p_{ia}}_i) -  \sum\limits_j \sum\limits_{b=1}^{l_j} CZ_{\tau} (\beta^{q_{jb}}_j).
\end{eqnarray*}
 By \cite{H1, H2} or  Section 3.4 of \cite{H4}, if $C$ is a simple holomorphic curve, then the ECH
index and Fredholm index satisfy the following \textbf{ECH index inequality}
\begin{equation} \label{eq20}
I(C) \ge \operatorname{ind}(C) +2\delta(C),
\end{equation}
where $\delta(C)$ is the count of singularities of $C$ with positive integer weights.
  Moreover,  the equality holds  only if $C$ is embedded and $C$ satisfies the \textbf{ECH partition conditions}. The ECH partition conditions mean that the partition of $C$ at
each $\alpha_i$ satisfies a certain special form. The general definition of such a form is quite
complicated. Here we only present the examples that will be considered in our proof. Suppose that $C$ has no negative ends and has positive ends at covers of an    Reeb orbit $\gamma$ with total multiplicities $m$. If $C$ satisfies the  {ECH partition conditions}, then the partition at $\gamma$ is
 \begin{itemize}
   \item
   $(1, ..., 1)$ if $\gamma$ is $L$-positive elliptic and $m$ satisfies $\mathcal{A}_{\lambda}(\gamma^m) <L$ ;
   \item
   $(m)$ if $\gamma$ is $L$-negative elliptic and $m$ satisfies $\mathcal{A}_{\lambda}(\gamma^m) <L$;
   \item
   $(1, ..., 1)$ $\gamma$ is positive hyperbolic.
 \end{itemize}

\paragraph{ECH group} Fix a class $\Gamma \in H_1(Y, \mathbb{Z})$. The chain group $ECC(Y, \lambda, \Gamma)$ is a free module generated by the ECH generators with homology class $\Gamma$.  Fix a generic  $J \in \mathcal{J}(Y, \lambda)$. The differential  is defined by
\begin{equation*}
\langle \partial \alpha, \beta \rangle: =\sum_{Z \in H_2(Y, \alpha, \beta), I(Z) =1} \# \left(\mathcal{M}^{J}(\alpha, \beta, Z)/ \mathbb{R}\right).
\end{equation*}
%Therefore, it induces a filtration on $ECH$.  More precisely,
Then $ECH(Y, \lambda, \Gamma)$ is the homology of this chain complex  $(ECC(Y, \lambda, \Gamma), \partial)$.

Given $L >0$, define $ECC^L(Y, \lambda, \Gamma)$ to be a submodule generated by the ECH generators with $\mathcal{A}_{\lambda}<L$. Note that the differential $\partial$ decreases the action. Therefore,   $ECC^L(Y, \lambda, \Gamma)$ is a subcomplex and its homology is well defined, denoted by $ECH^L(Y, \lambda, \Gamma)$. The group  $ECH^L(Y, \lambda, \Gamma)$ is called \textbf{filtered ECH}.  The inclusion induces a homomorphism
\begin{equation*}
i_L : ECH^L(Y, \lambda, \Gamma) \to ECH(Y, \lambda, \Gamma).
\end{equation*}

% Our goal is to compute the ECH cobordism maps $$ECH(E, \Omega): ECH(M, \lambda) \to \mathbb{Z}_2$$ induced by this filling.
\paragraph{ECH spectrum}
ECH also equips with a  homomorphism  $$U: ECH(Y, \lambda, \Gamma) \to ECH(Y, \lambda, \Gamma),$$ called the \textbf{U map}. It is defined by counting $I=2$ holomorphic currents passing through a fixed point (see (\ref{eq22})).

Suppose that $c_1(\xi)$ is torsion. Then $I(\alpha, \beta, Z)$ is independent of the choice of $Z$ for
null-homologous orbit sets $\alpha $ and $ \beta$. Hence, we write it as  $I(\alpha, \beta)$. For any orbit set $\alpha$ with $[\alpha]=0$, define its \textbf{grading} by
\begin{equation} \label{eq10}
\operatorname{gr}(\alpha):= I(\alpha, \emptyset).
\end{equation}
This gives a $\mathbb{Z}$ grading on $ECH(Y, \lambda, 0)$. Note that  the  $U$ map is a degree $-2$ map with respect to  this grading.
\begin{remark}
In the case that $Y$ is a prequantization bundle, by Lemma 3.11 of \cite{NW}, $c_1(\xi)$ is torsion. Therefore,     the grading (\ref{eq10}) on $ECH_*(Y, \lambda, 0)$ is well defined.
\end{remark}

There is  a canonical element $[\emptyset] \in  ECH(Y, \lambda, 0)$ which is represented by the empty orbit set.  The class $[\emptyset] $  is  called  the \textbf{contact invariant.} We remark that $[\emptyset] \ne 0$ if $(Y, \lambda)$ admits a symplectic filling.  So the  contact invariant of the prequantization bundle $\pi: Y \to \Sigma$ is nonzero.

Assume that $\lambda$ is nondegenerate. For $k \in \mathbb{Z}_{\ge 1}$, the \textbf{$k$-th value on ECH spectrum}
is defined by
\begin{equation*}
c_k(Y, \lambda) : = \inf\{L \in \mathbb{R}\vert     \exists   \sigma  \in ECH^L(Y, \lambda, 0) \mbox{ such that } U^k   (\sigma)= [\emptyset]\}.
\end{equation*}
If $\lambda$ is degenerate, define $c_k(Y, \lambda) : =  \lim_{n \to \infty} c_k(Y, f_n\lambda), $ where $f_n:Y\to \mathbb{R}_{>0}$  is a sequence of  smooth  functions such that $f_n\lambda$ is nondegenerate and $f_n$  converges to $1$ in $C^0$ topology.

\paragraph{Cobordism maps}
 A (strong) symplectic cobordism from $(Y_+, \lambda_+)$ to $(Y_-, \lambda_-)$
is a symplectic manifold $(X,\Omega_X) $ such that $\partial X =Y_+ \bigsqcup(-Y_-)$ and $\Omega_X(V_{\pm}, \cdot) =\lambda_{\pm}$ along $Y_{\pm}$, where $V_{\pm}$ are Liouville vector fields defined near $Y_{\pm}$. Using $V_{\pm}$, we can identify collar neighborhoods of $Y_{\pm}$ symplectically as $((-\delta, 0] \times Y_+, d(e^s \lambda_+) )$ and  $([0,\delta) \times Y_-, d(e^s \lambda_-) )$. We glue $(X,\Omega_X)$ with cylindrical ends $([0, \infty) \times Y_+, d(e^s \lambda_+))$ and
$(-\infty, 0]\times Y_-, d(e^s\lambda_-))$ along $Y_{\pm}$ to obtain a symplectic manifold $(\widehat{X},\Omega_X)$, which we call
the symplectic completion $(X,\Omega_X)$.  The ECH index, Fredholm index and holomorphic currents can be defined similarily in the cobordism setting (see \cite{H2}). Furthermore, the ECH index inequality (\ref{eq20}) still holds.

Fix a relative class $A \in H_2(X, \partial X, \mathbb{Z})$ such that $\partial A =\Gamma_+ -\Gamma_-$. Suppose that $\lambda_{\pm}$ are $\max\{L,  L+\rho(A)\}$-nondegenerate.  Here $\rho: H_2(X,\partial X, \mathbb{Z}) \to \mathbb{R}$ is a homomorphism defined by $\rho(A) : =\int_A \Omega_X -\int_{\Gamma_+} \lambda_+ + \int_{\Gamma_-} \lambda_-$. Hutchings and Taubes define a canonical homomorphism
\begin{equation*}
ECH^L(X, \Omega_X, A): ECH^L(Y_+, \lambda_+, \Gamma_+) \to ECH^{L+ \rho(A)}(Y_-, \lambda_-, \Gamma_-)
\end{equation*}
which is called a  \textbf{cobordism map} \cite{HT}.   If $\lambda_{\pm}$ are nondegenerate, then we can take $L \to \infty$ and get a cobordism map on the whole ECH groups
\begin{equation*}
ECH(X, \Omega_X, A): ECH(Y_+, \lambda_+, \Gamma_+) \to ECH(Y_-, \lambda_-, \Gamma_-).
\end{equation*}
The cobordism map $ECH^L(X, \Omega_X, A)$ is defined by counting the gauge classes of solutions to the Seiberg-Witten equations. We will not provide any details about the Seiberg-Witten theory here, and  we refer readers to  the book of   P. Kronheimer  and T. Mrowka \cite{KM}. Hutchings and Taubes show that $ECH^L(X, \Omega_X, A)$ satisfies the \textbf{holomorphic curve axioms} (see  Theorem 1.9 of  \cite{HT}). Roughly speaking, it means that if the cobordism map is nonvanishing, then  there exists a holomorphic current.   In some special cases, C. Gerig enhances the holomorphic curves axioms. He shows that there is a one-to-one  correspondence between the holomorphic currents and the gauge classes of solutions to the Seiberg-Witten equations \cite{CG, CG2}. In other words, the cobordism map is actually defined by counting holomorphic curves. We will show that this is the case in our situations.

\section{Computations of the  cobordism maps} \label{section3}
In this section, we prove the Theorem  \ref{thm2}.
\subsection{Perturbations}\label{section31}
Before we go ahead, we need to clarify the perturbations made on the contact form $\lambda$ and the symplectic form $\Omega$.
\paragraph{Morse-Bott perturbations}
Note that the contact form $\lambda$ is Morse-Bott. The Reeb orbits are iterations  of the fibers of $\pi: Y \to \Sigma$. Following   Farris, Nelson and Weiler's approach  \cite{Fa, NW},  we  perturb the contact form by   a perfect   Morse function $H: \Sigma \to \mathbb{R}$. More precisely,  define
\begin{equation*}
\lambda_{\varepsilon}: = (1+{\varepsilon}\pi^*H) \lambda,
\end{equation*}
where   $0 < \varepsilon \ll1  $ is a small fixed number.

Let $e_-$, $e_+$ and  $\{h_i\}_{i=1}^{2g}$ denote the fiber over the minimum, maximum and saddle  points of $H$ respectively. These are simple  Reeb orbits of $\lambda_{\varepsilon}$.    Moreover,  $e_{\pm}$ are elliptic orbits and $\{h_i\}_{i=1}^{2g}$ are positive hyperbolic orbits.

For any $0<\varepsilon \ll1$, there exists a constant $L_{\varepsilon}$ such that   $ \lambda_{\varepsilon}$ is $L_{\varepsilon}$-nondegenerate and the covers of  $e_{\pm}$, and $\{h_i\}_{i=1}^{2g}$ are  the only  Reeb orbits of $(Y, \lambda_{\varepsilon})$ with action less than  $L_{\varepsilon}$.  According to the computations  in Lemma 3.9 of \cite{NW}, $e_+$ is $L_{\varepsilon}$-positive and $e_-$ is $L_{\varepsilon}$-negative.
We remark that $L_{\varepsilon} \to \infty$ as $\varepsilon \to 0$.  For an orbit set $\alpha=e_{+}^{m_+}h_1^{m_1}\cdots h_{2g}^{m_{2g}} e_-^{m_-}$,  its action is given by
\begin{equation*}
\mathcal{A}_{\lambda_{\varepsilon}}(\alpha) =2M + \varepsilon \left(m_+ \pi^*H(e_+) + m_- \pi^*H(e_-)+  \sum_{i=1} m_i \pi^*H(h_i)  \right),
\end{equation*}
where $M =m_- + m_++ \sum_{i=1}^{2g} m_i.$

By Lemma 3.7 in \cite{NW}, we have  $H_1(Y, \mathbb{Z}) = \mathbb{Z}^{2g} \oplus \mathbb{Z}_{|e|}$.  The homology class of each fiber of $\pi: Y \to \Sigma$ is $1 \mod |e|$ in the  $\mathbb{Z}_{|e|}$ summand of $H_1(Y, \mathbb{Z})$.  Therefore,   $M =d|e| + \Gamma$ for some integer $d\ge 0 $ and $\Gamma =[\alpha] \in H_1(Y, \mathbb{Z}).$

\begin{comment}
Define a diffeomorphism $F: \mathbb{R}_{>0}\times Y \to \mathbb{R}_{>0} \times Y$ by
\begin{equation*}
F(t, x) = (t(1+\varepsilon(t)\pi^*h),x).
\end{equation*}
Then $F^*(r\lambda)=tf\lambda$ when $t > \frac{4}{5}$. Under the identification \ref{eq1}, we can extend $F$ as a self diffeomorphism of $E$.
\end{comment}

% $A=\frac{1}{2}(\bar{z} dz -  zd\bar{z} )$. We choose $H$ such that $H=\varepsilon |z|^2$ over $U$.

\paragraph{Symplectic completion of $(DE, \Omega)$} We regard $E$ as a symplectic completion of
$(DE, \Omega)$ by the following way. Define $r :=\frac{1}{2}(1 + \rho^2)$. Then we have a symplectomorphism
\begin{equation} \label{eq1}
  (E\setminus \Sigma, \Omega) \cong (\mathbb{R}_{r > \frac{1}{2}} \times Y, d(r\lambda)).
\end{equation}
Sometimes we identify the conical end $(\mathbb{R}_r  \times  Y, d(r\lambda))$ with the cylindrical end $(\mathbb{R}_s \times Y, d(e^s \lambda))$ via changing coordinate $r = e^s$.

Let $\varepsilon(r)$ be a nondecreasing cutoff function such that $\varepsilon(r)=\varepsilon \ll 1$ when  $r\ge  \frac{4}{5}$ and  $\varepsilon(r)=0$ when $r \le \frac{3}{4}$. Define $\lambda_{\varepsilon(r)}: =(1+ \varepsilon(r) \pi^*H) \lambda. $   Under the identification (\ref{eq1}), we define
\begin{equation} \label{eq15}
\Omega_{\varepsilon} =
\begin{cases}
\Omega, & r  <  \frac{3}{4}\\
d(r \lambda_{\varepsilon(r)}),  & r\ge \frac{3}{4}.
\end{cases}
\end{equation}
 Note that $\Omega_{\varepsilon} = d(e^s \lambda_{\varepsilon})$  over the cylindrical end.

\paragraph{Almost complex structures in $E$} An $\Omega_{\varepsilon}$-compatible almost complex structure $J$ is \textbf{cobordism admissible} if  $J=J_+$ for some $J_+ \in \mathcal{J}(Y, \lambda_{\varepsilon})$ over the cylindrical end.   We choose $J$ such that
\begin{enumerate} [label=\textbf{ J.\arabic*}]
   \item \label{J1}
In the neighbourhood $U_z \times \mathbb{C}_w$ of $e_-$.     $J(r \partial_r) = f(r)\partial_{\theta}$ and $J(\partial_{z}) = i \partial_z$ along  the fiber $\{0\} \times \mathbb{C}$, where $f(r)$ is a positive function such that $f(r) = 1$ when $r \ge e^{-\epsilon}$ and $f(r) = \frac{r}{2r-1}$ when $r $ is close to   $\frac{1}{2}$. The latter assumption on $f$ is equivalent to $J(\rho \partial_{\rho}) =\partial_{\theta}$. This implies that $J$ is well defined on the whole $\{0\} \times \mathbb{C}$.
   \item \label{J2}
   The zero section $\Sigma$ is $J$-holomorphic.
 \end{enumerate}
The choice of $J$ satisfying  \ref{J1} and  \ref{J2} is always feasible.
The advantage of \ref{J2} is
that $\mathcal{C}\cdot \Sigma \ge 0$  for any holomorphic current $\mathcal{C}$ without closed components. We use this
property in Lemma \ref{lem9} and Lemma \ref{lem4}.

 By \ref{J1},  $C_{e_-} =\{0\} \times \mathbb{C}$ is  holomorphic because $J(TC_{e_-} ) \subset TC_{e_-} .$   Moreover, $C_{e_-}$ is asymptotic to $e_-$ under the identification (\ref{eq1}). We remark that $C_{e_-}$ is Fredholm regular for any $J$. This follows directly from  C. Wendl's automatic transversality theorem \cite{Wen} and the index computation in Lemma \ref{lem2}.  The fiber $C_{e_-}$ plays the same  role as the horizontal section in \cite{GHC1}.

An almost complex structure  $J$ is called \textbf{generic} if all simple holomorphic curves are Fredholm regular except for the closed holomorphic curves. We assume that $J$ is generic, unless   stated otherwise.

\subsection{Moduli space of holomorphic currents}
\paragraph{Computing the ECH index} The first task of studying the holomorphic currents is to compute their index. The computations here follow the similar argument  in \cite{GHC1}, where the author computes ECH index  of the relative homology classes in an elementary Lefschetz fibration. This approach coincides with  Nelson and Weiler's methods (see Remark \ref{remark2}).

%Given an orbit set $\alpha$, let $H_2(DE, \alpha)$ denote the set of relative homology classes.  It is an affine space over $H_2(DE, \mathbb{Z})$.

Let $\alpha=e_{+}^{m_+}h_1^{m_1}  \cdots  h_{2g}^{m_{2g}} e_-^{m_-}$ be an orbit set.   Let $M:=m_+ + m_1 + \cdots + m_{2g}+ m_-$.    Let $H_2(DE, \alpha)$ denote the set of relative homology classes.   It is an affine space over $H_2(DE, \mathbb{Z})$.   Let $C_{e_{\pm}}$ and $C_{h_i}$ denote the fibers over the critical points corresponding to $e_{\pm}$ and $h_i$. For $\alpha$, we define a relative homology  class $Z_{\alpha}$ represented by
\begin{equation*}
m_+ C_{e_+} + \sum_{i=1}^{2g}m_i C_{h_i}+ m_- C_{e_-}.
\end{equation*}
Since $H_2(DE, \mathbb{Z}) \cong H_2(\Sigma, \mathbb{Z})$ is generated by $[\Sigma]$, any other relative homology class can be written as $Z_{\alpha} + d [\Sigma]$ for some integer $d$.

Fix a relative homology class $A \in H_2(DE, Y, \mathbb{Z})$.  Consider a subset $\{Z \in H_2(DE, \alpha): [Z]  =A \}$.  Note that it is an affine space of $\ker j_*$, where $j_*$ is the map in the following exact sequence
 \begin{equation*}
\begin{split}
...&\to  H_2(Y, \mathbb{Z}) \xrightarrow{ i_*}H_2(DE,  \mathbb{Z}) \xrightarrow{ j_*} H_2(DE, Y, \mathbb{Z})  \xrightarrow{ \partial_*} H_1(Y, \mathbb{Z}) \to....
\end{split}
 \end{equation*}
By Lemma 3.7 in \cite{NW},  $H_2(Y, \mathbb{Z}) $ is spanned  by the classes  represented by $\pi^{-1}(\eta)$, where $\eta$ is a simple closed curve in $\Sigma$.  Note that  $\pi^{-1}(\eta)$   is the boundary of  the surface $\pi_E^{-1}(\eta) \cap DE$ in $DE$.   Then,  $i_*$ is  the  zero map and $\ker j_*=0$.  Therefore, a relative class  $A \in H_2(DE, Y , \mathbb{Z})$ and an orbit set $\alpha$ determine  a unique relative homology class $Z_{\alpha, A} \in H_2(DE, \alpha)$.

For a relative homology class  $Z \in H_2(DE, \alpha)$, we obtain the following formula for the ECH index.
\begin{lemma}\label{lem1}
Given an orbit set  $\alpha=e_{+}^{m_+}h_1^{m_1}\cdots h_{2g}^{m_{2g}} e_-^{m_-}$, the ECH index of a relative class $Z_{\alpha}+ d[\Sigma] \in H_2(DE, \alpha)$ is
\begin{equation*}
I(Z_{\alpha}+d[\Sigma])=M+m_+-m_- + 2dM +d^2 e +de + d\chi(\Sigma).
\end{equation*}
\end{lemma}
\begin{proof}
Let $p$ be a critical point of $H$, $\gamma_p=\pi^{-1}(p)$ and $C_p=\pi^{-1}_E(p)$. We fix a constant trivialization as in \cite{NW}. More precisely, a trivialization of $T_p\Sigma$ can be lifted to a trivialization of $\xi\vert_{\gamma_p(t)}$.  Under this trivialization, by Lemma 3.9 of \cite{NW}, we have
\begin{equation*}
CZ_{\tau}(\gamma_p^k) = \operatorname{ind}_p H-1,
\end{equation*}
where $\operatorname{ind}_p H$ is the Morse index at $p$.

Regard the connection $A_{\nabla}$ as a map $A_{\nabla}: TE \to  E$. Then it  induces a splitting $TE=TE^{hor}\oplus TE^{vert}$, where $TE^{hor} : = \ker A_{\nabla}$ and $TE^{vert} := \ker d\pi_E$.
The trivialization of $T_p\Sigma$ also can be lifted to $T^{hor} E \vert_{C_p}$. In particular, $c_{\tau}(T^{hor} E \vert_{C_p})=0$. Note that  $T^{hor} E \vert_{C_p}$ can be identified with the normal bundle $N_{C_p}$ of $C_p$. By the definition of $Q_{\tau}$, we obtain  $Q_{\tau}(C_p)=c_{\tau}(N_{C_p})=0$.  The  section  $r \partial_r -i\partial_{\theta}$   can be extended to  a section  $\psi$ of $T^{vert} E \vert_{C_p} =TC_p$. We choose such an extension  such that  $\psi= \rho \partial_{\rho}- i\partial_{\theta}$ near $r =\frac{1}{2}$  and   is non-vanishing except at the zero section.  Hence,  the relative Chern number is $$c_{\tau}( T^{vert} E \vert_{C_p})= \#\psi^{-1}(0) = \#( \rho \partial_{\rho}- i\partial_{\theta})^{-1}(0)=1.$$
Therefore,  we have
\begin{equation*}
I(C_p)=c_{\tau}(TE\vert_{C_p}) + Q_{\tau} (C_p) + CZ_{\tau}(\gamma_p)=\operatorname{ind}_p H.
\end{equation*}
In other words, $I(C_{e_+})=2$, $I(C_{h_i})=1$ and  $I(C_{e_-})=0$. For $p\ne q$, the fibers $C_p$ and $C_q$
are obviously disjoint. Hence, $Q_{\tau}(C_p, C_q)=0$. Therefore,
the formula for $I(Z_{\alpha})$ follows from the facts that the relative Chern number is additivity  and $Q_{\tau}$ is quadratic.

By the definition of the ECH index, we have
\begin{equation}\label{eq2}
I(Z_{\alpha}+ d[\Sigma]) =I(Z_{\alpha}) + 2d Z_{\alpha} \cdot \Sigma +I(d\Sigma)= M-m_- + m_+ + 2dM +I(d\Sigma).
\end{equation}
Under the identification $T^{hor}E=\pi_E^* T\Sigma$ and $T^{vert}E=\pi_E^*E$, we have
\begin{equation*}
\langle c_1(TE), \Sigma \rangle= \langle c_1(T^{hor}E), \Sigma \rangle + \langle c_1(T^{vert}E), \Sigma \rangle=\chi(\Sigma) + e.
\end{equation*}
The self-intersection number of $\Sigma$ is
\begin{equation*}
\Sigma\cdot \Sigma =\#s^{-1}(0)=\langle c_1(E), \Sigma \rangle=e,
\end{equation*}
where $s$ is a generic section of $E$. In sum, we have
\begin{equation}\label{eq3}
I(d\Sigma)=d\chi(\Sigma) + de + d^2e.
\end{equation}
The statement of the lemma follows from  (\ref{eq2}) and (\ref{eq3}).
\end{proof}
\begin{remark} \label{remark2}
Let $\alpha =e_+^{m_+} h_1^{m_1},...,h_{2g}^{m_{2g}} e_-^{m_-}$ and $\beta =e_+^{n_+} h_1^{n_1},...,h_{2g}^{n_{2g}} e_-^{n_-}$. In Proposition  3.5 of \cite{NW},  Nelson and Weiler show that
\begin{equation} \label{eq8}
I(\alpha, \beta)= d\chi(\Sigma) + d^2|e| +2dN +m_+ - m_- + n_--n_+,
\end{equation}
where $M$ and $N$ are total multiplicities of $\alpha$ and  $\beta$ respectively, and $d = (M - N)/|e|$.
 Our computations here is equivalent to (\ref{eq8})     in the following sense.

  One can use Lemma \ref{lem1} and the additivity of ECH index to recover the index formula $I(\alpha, \beta)$ as follows.  Let $Z_{\alpha} + d_{\alpha} [\Sigma] \in H_2(DE, \alpha)$, $Z_{\beta} + d_{\beta} [\Sigma] \in H_2(DE, \beta)$  and $Z \in H_2(Y, \alpha, \beta)$ be relative homology classes such that  $Z_{\alpha} + d_{\alpha} [\Sigma]  =Z \#(Z_{\beta} + d_{\beta} [\Sigma])$.  By a similar  argument  as in Proposition  4.4 of \cite{NW} (see (\ref{eq12}) for the details in the current setting),
one can show that $d_{\alpha} -d_{\beta} =(M-N)/|e|$. Then, by direct computations, we have
\begin{equation}\label{eq14}
I(\alpha, \beta) =I(Z) = I(Z_{\alpha} + d_{\alpha} [\Sigma] ) -I(Z_{\beta} + d_{\beta} [\Sigma]).
\end{equation}
Following from Lemma \ref{lem1} and (\ref{eq14}), we get the formula (\ref{eq8}) for $I(\alpha, \beta)$.

Conversely, we also can deduce Lemma \ref{lem1} easily by applying  Nelson and  Weiler's index computations as follows.  Suppose that   $[\alpha]=0$ and $M=k|e|$.  Let $Z_0=Z_{\alpha} +k[\Sigma]$. Note that $Z_0 \cdot [\Sigma]=0$. We can represent $Z_0$ by a surface $S_0$ which is disjoint from $\Sigma$. Under the identification (\ref{eq1}), $S_0$ represents a relative homology class  in $H_2(Y, \alpha, \emptyset)$. In particular, we have $I(Z_0)=I(S_0)=I(\alpha, \emptyset)$.  Then one can apply the index ambiguity formula (\ref{eq2}) to drive the result for a general $Z=Z_0 +d[\Sigma]$.
\end{remark}

In the next lemma, we use the computations in Lemma \ref{lem1} to drive a formula for the Fredholm index.
\begin{lemma} \label{lem2}
Let $C$ be a holomorphic curve with relative class $[C]=Z_{\alpha}+ d[\Sigma] $. Then the  Fredholm index  of $C$ is
\begin{equation}
\operatorname{ind}(C)=2g(C)-2 + h(C) + 2e_+(C)+ 2M+2d\chi(\Sigma) +2de,
\end{equation}
where $h(C)$ is the number of ends at hyperbolic orbits and $e_+(C)$ is the number of ends at covers of $
e_+$.
\end{lemma}
\begin{proof}
The proof follows directly from definition and the calculations  in Lemma \ref{lem1}.
\end{proof}

\paragraph{$J$-Holomorphic currents without closed components}We first  study  the case  that the $J$-holomorphic current $\mathcal{C}$ contains no closed components. Also, we assume that the holomorphic curves are asymptotic to orbit sets with $\mathcal{A}_{\lambda_{\varepsilon}} < L_{\varepsilon}$,  unless stated otherwise.

Fix a relative homology class $Z \in H_2(DE, \alpha)$. Let $\mathcal{M}^J (\alpha, Z)$ denote the set  of $J$-holomorphic currents in $E$. We denote the moduli space of \textbf{broken $J$-holomorphic currents} in $E$ by $\overline{\mathcal{M}^J (\alpha, Z)}$  (Definition 1.6 of \cite{HT}).  An element in  $\overline{\mathcal{M}^J (\alpha, Z)}$ is a chain of holomorphic currents $\{\mathcal{C}^0, ..., \mathcal{C}^{N}\}$ such that
\begin{enumerate}
  \item
  $\mathcal{C}^0 \in \mathcal{M}^J(\alpha_0, Z_0)$ and   $\mathcal{C}^i \in \mathcal{M}^J(\alpha_i, \beta_i)$ for $i \ge 1$;
  \item
  $\alpha_i =\beta_{i+1}$ for $1\le i\le N-1$, and $\alpha^N =\alpha$;
  \item
 $ Z= \sum_{i=0}^N[\mathcal{C}^i] $.
\end{enumerate}
By Lemma 6.8 of \cite{HT}, a sequence of  $J$-holomorphic currents in ${\mathcal{M}^J (\alpha, Z)}$ converges to
a broken $J$-holomorphic current in  $\overline{\mathcal{M}^J (\alpha, Z)}$ in the current sense (see Section 9 of \cite{H1} for details).

For a $J$-holomorphic curve in $E$, we define its degree which is an analog of Definition
4.1 in \cite{NW}.  Let  $\mathcal{C} \in \mathcal{M}^J(\alpha, Z)$ be  a $J$-holomorphic current represented by a $J$-holomorphic map $ u: \dot{F} \to E$, where $\dot{F} =F \setminus P$, $F$ is a closed  Riemann surface (possibly disconnected) and $P$ is the set of punctures.  Since $\pi_E \circ u$ maps the punctures to the critical points of $H$, we  extend $\pi_E \circ u$  to a map $\pi_E \circ u: F \to \Sigma$. Then, we have a well-defined degree $\operatorname{deg}(\pi_E \circ u)$. Define $\operatorname{deg}(\mathcal{C}) :=\operatorname{deg}(\pi_E \circ u)$. It is called the  \textbf{degree}  of $\mathcal{C}$.  Alternatively, we can define $\operatorname{deg}(\mathcal{C})$ to be the unique integer $d$ such that $[\mathcal{C}] =Z_{\alpha} +d [\Sigma]$. To understand why this is valid, consider the following computation: $\int_{\mathcal{C} } \pi_E^*\omega_{\Sigma} = \int_{\pi_E(\mathcal{C})} \omega_{\Sigma} = \operatorname{deg}(\mathcal{C})|e|$. On the other hand,   $\int_{\mathcal{C} } \pi_E^*\omega_{\Sigma} = \int_{Z_{\alpha} + d[\Sigma]} \pi_E^*\omega_{\Sigma} = d|e|.$ Thus, we conclude that $d = \operatorname{deg}(\mathcal{C})$. The degree of a broken holomorphic curve is
defined in the same way.

\begin{lemma} \label{lem11}
 For a generic almost complex structure $J$, let $\mathcal{C}$ be a $J$-holomorphic current without closed components.  Then,   we have $\operatorname{deg}(\mathcal{C}) \ge 0$.
\end{lemma}
\begin{proof}
Write $\mathcal{C}=\sum_a d_a C_a.$  Since $\operatorname{deg}(\mathcal{C}) = \sum_ad_a \operatorname{deg} (C_a)$, it suffices to prove the conclusion for an irreducible simple holomorphic curve $C$ with at least one end.

Let $d=\operatorname{deg}(C)$. Assume that $d \le -1$. Then,   Lemma \ref{lem1} implies that
\begin{equation*}
 I(C) \le |d|(1-|d|)|e| + d\chi(\Sigma) \le d \chi(\Sigma).
\end{equation*}
 If $\chi(\Sigma)=2$, then  $I(C) \le -2$.
If $\chi(\Sigma)\le 0$, then  $\operatorname{ind} (C) \ge 2d\chi(\Sigma) + 2de> d\chi(\Sigma)$ by Lemma \ref{lem2}.

In both cases, they violate the ECH inequality $I(C) \ge \operatorname{ind}(C) \ge 0$.
\end{proof}

From now on, we shall assume that $\operatorname{deg}(C) \ge 0$ due to the above lemma.

To deal with the $J$-holomorphic currents with multiply covered components, we
need the following self-intersection number that appears in the ECH index inequality.
\begin{definition}[Definition 4.7 \cite{H5}]
For two  simple $J$-holomorphic curves $C, C'$ which are asymptotic positively to an orbit set $\alpha$ with  action less than $L$, define an integer $C \star C'$ as follows.
\begin{itemize}
\item
If $C$ and $C'$ are distinct, then  $C \star C'$
is the algebraic count of intersections of $C$ and $C'$.  By intersection positivity, we have $C \star C' \ge 0$. The equality holds if and only if $C$ and $C'$ are disjont.
\item
If $C$ and $C'$ are the same curve, then  define
\begin{equation}\label{eq16}
2C\star C=2g(C)-2+ h(C) + \operatorname{ind}(C) + 2e_L(C) + 4 \delta(C) ,
\end{equation}
where $e_L(C)$ is  the total
multiplicity of all elliptic orbits in $\alpha$ that are $L$-negative,  and $\delta (C)$ is the count of singularities of $C$ with positive integer weights.  $\delta(C) \ge 0$ and equality holds if and only if $C$ is embedded.
\end{itemize}
\end{definition}

Let $\mathcal{C}=\sum_{a} d_a C_a$ and $\mathcal{C}'=\sum_b d_b' C_b'$ be two $J$-holomorphic currents.   By Proposition 4.8 of \cite{H5}, we have
\begin{equation} \label{eq11}
I(\mathcal{C}+ \mathcal{C}') \ge I(\mathcal{C}) +I(\mathcal{C}') +2 \mathcal{C} \star \mathcal{C'},
\end{equation}
where $\mathcal{C} \star \mathcal{C'} = \sum_a \sum_b d_a d_b' C_a \star C_b'. $

If the self-intersection number (\ref{eq16}) is nonnegative, then the ECH index inequality  (\ref{eq20}) and (\ref{eq11}) imply that the ECH index is nonnegative.  However, this is not always the case. See Section 5.5 of \cite{H4} for the explanations.  In our setting, using the index computations in Lemma \ref{lem1} and Lemma \ref{lem2}, we show that if $C \star C<0$, then its ECH index is at least $2$.  As a result, the ECH indices  of the holomorphic currents in $E$ are nonnegative.

\begin{lemma} \label{lem9}
Let $C$ be an irreducible simple  holomorphic curve with at least one end. If $C \star C < 0$,
%\begin{equation*}
%2C\star C=2g(C)-2+ h(C) + ind C + 2e_L(C) + 4 \delta(C) < 0,
%\end{equation*}
%where $e_L(C)$ be the total
%multiplicity of all elliptic orbits in $\alpha$ that are $L$-negative,
then $I(kC)\ge 2$ for any $k \ge 1$. In particular, if $\mathcal{C}$ is a $J$-holomorphic
current without closed component, then $I(\mathcal{C}) \ge 0$.
\end{lemma}
\begin{proof}
Assume that $C \star C <0$. By Lemma \ref{lem2}, we know that $h(C)+ \operatorname{ind}(C)$ is a nonnegative even integer. By the definition (\ref{eq16}), we have  $$g(C)=h(C) =\operatorname{ind}(C) =e_L(C) =\delta(C)=0.$$ The condition $h(C)=e_L(C)=0$ forces $M=m_+$.

By Lemma \ref{lem2}, $\operatorname{ind}(C)=0$ implies that
\begin{equation} \label{eq4}
M=d|e|-d\chi(\Sigma)+ 1 -e_+(C).
\end{equation}
Write the relative homology class of $C$ as $Z_{{e_{+}^M}} + d[\Sigma]$, then $[kC]= Z_{e_+^{kM}} + dk[\Sigma]$. Note that $C \cdot \Sigma=M + de \ge 0$ by intersection positivity; then $M \ge d|e|$.  If $\chi(\Sigma)>0$, then
\begin{equation*}
\begin{split}
I(kC)&= 2kM +d^2k^2e + dke+dk\chi(\Sigma) +2dk^2M\\
&\ge 2kM + dk \chi(\Sigma) + k^2d^2|e| -dk|e| \ge 2kM \ge 2.
\end{split}
\end{equation*}

If $\chi(\Sigma) \le 0$, then  by Lemma \ref{lem2} and Equation (\ref{eq4}), we have
\begin{equation*}
\begin{split}
I(kC)&= 2kM +d^2k^2e + dke+dk\chi(\Sigma) +2dk^2M\\
&=kM + k(1+d|e|-d\chi(\Sigma)-e_+(C)) +d^2k^2 e + dk e+ dk\chi(\Sigma) \\
&+ dk^2M + dk^2(d|e|-d\chi(\Sigma) +1 -e_+(C))\\
&=k(M-e_+(C)) + k + dk^2(M-e_+(C)) -d^2k^2 \chi(\Sigma) + dk^2 \ge k \ge 1.
\end{split}
\end{equation*}
Note that $I(kC)$ is even. Hence, we get $I(kC) \ge 2$.

Write $\mathcal{C} =\sum_a d_aC_a + \sum_{b} d_b' C_b'$ such that $C_a \star C_a\ge 0$  and $C_b'\star C_b'<0$. By Inequality (\ref{eq11}), we have
\begin{equation}\label{eq5}
\begin{split}
I(\mathcal{C}) \ge& \sum_a I(d_aC_a) +  \sum_b I(d_b'C_b') + 2\sum_{a,b} d_a d_b' C_a \star C_b' + 2\sum_{a\ne a'} d_a d_{a'} C_a \star C_{a'} + 2\sum_{b\ne b'} d_b' d_{b'}' C_b' \star C_{b'}'\\
\ge & \sum_a d_aI(C_a) + \sum_a d_a(d_a-1)C_a \star C_a +   \sum_b I(d_b'C_b') + 2\sum_{a,b} d_a d_b' C_a \star C_b' \\
+& 2\sum_{a\ne a'} d_a d_{a'} C_a \star C_{a'} + 2\sum_{b\ne b'} d_b' d_{b'}' C_b' \star C_{b'}' .\\
\end{split}
\end{equation}
By ECH index inequality (\ref{eq20}),  $\sum_a d_aI(C_a) \ge 0$.   The intersection positivity implies that the last three terms are nonnegative.  $C_a \star C_a \ge 0$ by assumption and $ I(d_b'C_b')  \ge 2$. In sum, $I(\mathcal{C}) \ge 0$.
\end{proof}

A simple $J$-holomorphic curve $C$ is called a \textbf{special $J$-holomorphic plane} if it has $I(C) =\operatorname{ind}(C) =0$, and is   an embedded plane whose positive end  is
asymptotic to $e_-$ with multiplicity $1$. This is a counterpart of the Definition 3.15 in \cite{CG}.
\begin{lemma} \label{lem4}
Assume that $C$ is not closed.  If $I(C)=\operatorname{ind}(C) =C\star C =0$, then $C$ is a special $J$-holomorphic plane.
\end{lemma}
\begin{proof}
By Definition (\ref{eq16}),    $C \star C=0$ forces $\delta(C)=0$, i.e, $C$ is embedded.  It is easy to check that  $C $ satisfies one of the following  properties:
 \begin{enumerate}
  \item
 $h(C)=e_L(C)=0$ and $g(C)=1$;
 \item
 $h(C)=2$ and $g(C)=e_L(C)=0$;
 \item
 $h(C)=g(C)=0$ and $e_L(C)=1$.
 \end{enumerate}
%In the first case, all the ends of $C_a$ are asymptotic to covers of $\gamma_+$. The total multiplicity of the covers is denoted by  $M_a$.
Let $d =\operatorname{deg} (C)$. By Lemma \ref{lem1} and  Lemma  \ref{lem2}, we have
\begin{equation} \label{eq6}
\begin{split}
2I(C)-\operatorname{ind}(C)=2m_+ - 2m_-  +4dM -2d^2|e| +2 -2g(C)-h(C) -2e_+(C).
\end{split}
\end{equation}
Since $I(C)=\operatorname{ind}(C)$, the ECH partition condition implies that $e_+(C)=m_+$. Also, $e_L(C) =0$ is equivalent to $m_-=0$.   In the first two cases,  we have
\begin{equation*}
\begin{split}
0=4dM -2d^2|e| \ge 2d M \ge 0.
\end{split}
\end{equation*}
The last step comes from the positivity intersection of holomorphic curves $C \cdot \Sigma =M-d |e| \ge 0 $. Hence, we have either $M=0$ or $d=0$. If $d=0$, then the formula in Lemma  \ref{lem1} still implies that  $M=0$. We get contradiction since we have assumed that $C$ is not closed.

In the last case, $m_-=1$.  By Equation (\ref{eq6}), then we still get $d=0$. The formula in Lemma \ref{lem1} and $ I(C)=0$ imply that $m_+=0$. Hence, $C$ is a $J$-holomorphic plane with one end at $e_-$, i.e., it is a special $J$-holomorphic plane.
\end{proof}

\begin{lemma} \label{lem5}
 Let $\mathcal{C}\in \mathcal{M}^J(\alpha, Z)$ be a $J$-holomorphic current with $I(\mathcal{C}) =i$, where $i=0$ or $1$. If $i=1$, we also  assume that $\alpha$ is an ECH generator.   Suppose that $\mathcal{C}$ has no closed component. Then  $\mathcal{C} =\mathcal{C}_{emb} \cup \mathcal{C}_{spec}$, where $\mathcal{C}_{emb}$ is embedded with $I(\mathcal{C}_{emb})= \operatorname{ind} (\mathcal{C}_{emb})=i$ and $\mathcal{C}_{spec}$ consists of special $J$-holomorphic planes.
\end{lemma}
\begin{proof}
 Write $\mathcal{C} =\sum_a d_aC_a + \sum_{b} d_b' C_b'$  as in the proof of  Lemma \ref{lem9},  with the properties that $C_a \star C_a \ge 0$ and $C_{b}' \star C_{b}'<0$. By Lemma \ref{lem9} and Inequality (\ref{eq5}), we must have $d_b'=0$  because $I(d_b'C_b')\ge 2$.

In the case that $I(\mathcal{C})=0$,   we have $I(C_a)=0$ for any $a$. Also, $d_a=1$ unless $C_a \star C_a=0$. The ECH index inequality  implies that $\operatorname{ind} (C_a)=0$ and $\delta(C_a)=0$ as well. If $d_a>1$, then  $C_a \star  C_a=0$.  By Lemma \ref{lem4}, $C_a$ is a special $J$-holomorphic plane.

In the case that $I(\mathcal{C})=1$, then we have $I(C_a) \le 1$. If $I(C_a)=0$ for all $a$, the ECH index inequality and Lemma  \ref{lem2} implies that $C_a$ has even ends at hyperbolic orbits. Since $\alpha$ is an  ECH generator, we know that $\alpha$ contains even distinct  simple hyperbolic orbits. By Lemma \ref{lem1}, $I(\mathcal{C})= 0\mod 2$, we get a contradiction.  Therefore, there exists $C_{a_0}$ with $I(C_{a_0})=\operatorname{ind} (C_{a_0}) =1$. The Inequality (\ref{eq5}) implies that such $a_0$ is unique and $d_{a_0}=1$. For any other $a$, we also have  $I(C_a)=\operatorname{ind} (C_a)= \delta(C_a)=0$. Moreover,  $d_a=1$ unless $C_a \star C_a=0$.

In both cases, $\mathcal{C}$ is a union of embedded curves and covers of special $J$-holomorphic planes.

% Then $ind C_a=0$ and ECH partition condition imply that
%\begin{equation*}
%2M_a= -d_a\chi(B) + d_a |e|.
%\end{equation*}
%If $d_a=0$, then $C_a$ is closed. We can assume that $d_a \ge 1$. By Lemma \ref{lem1}, we have
%\begin{equation*}
%\begin{split}
%I(C_a)&=2M_a+2d_a M_a -d_a^2|e|-d_a|e|+d_a \chi(B)\\
%&=2d_a M_a -d_a^2|e| \ge d_a M_a \ge 1.
%\end{split}
%\end{equation*}
%The last step comes from the positivity intersection of holomorphic curves $C_a \cdot B =M_a-d_a |e| \ge 0 $. But this contradicts with $I(C_a)=0$.
\end{proof}

\paragraph{Closed holomorphic curves} Now we begin to consider the holomorphic currents
that contain closed $J$-holomorphic curves. We first need to figure out what kind of
closed $J$-holomorphic curves could exist in $E$.
\begin{lemma} \label{lem7}
The zero section $\Sigma$ is the unique simple  closed $J$-holomorphic curve in $E$.
\end{lemma}
\begin{proof}
Suppose we have a  simple closed holomorphic curve $C$ which is different from  $\Sigma$. Since $H_2(DE, \mathbb{Z})$ is generated by $[\Sigma]$, we must have $[C] =k [\Sigma]$. By energy reason, we have  $k \ge 1$. However, $C\cdot \Sigma =k[\Sigma] \cdot [\Sigma] =ke <0$, contradicts with the intersection positivity of holomorphic curves.
\end{proof}

\begin{lemma} \label{lem6}
%Fix a relative class $A \in H_2(E, Y, \mathbb{Z})$ represented by $Z_{e_-^M}$.
Let $\mathcal{C} \in \overline{\mathcal{M}^J(\alpha, Z_{\alpha})}$ be a  broken $J$-holomorphic current . Then $\mathcal{C}$ contains no closed components. Moreover, each level of $\mathcal{C}$ has degree $0$.
\end{lemma}
\begin{proof}
By  Lemma \ref{lem7},  we can write $\mathcal{C} = \mathcal{C}_{\diamondsuit} + k \Sigma$, where $ \mathcal{C}_{\diamondsuit}  $ has no closed components and $k \ge 0$.

 Write $ \mathcal{C}_{\diamondsuit} = \{\mathcal{C}^0,...,\mathcal{C}^N\}$, where $\mathcal{C}^0 \in \mathcal{M}^J(\alpha_0)$, $\mathcal{C}^i \in \mathcal{M}^J(\alpha_i, \alpha_{i-1})$ and $\alpha^N =\alpha.$ We
claim that the degree is additive, i.e.,  $\operatorname{deg}(\mathcal{C}_{\diamondsuit}) =\sum_{i=0}^N \operatorname{deg}(\mathcal{C}_{i}) $, where $\operatorname{deg}(\mathcal{C}_{i}) (i \ge 1)$  is the degree defined in Definition 4.1 of \cite{NW}. Proposition 4.4 of \cite{NW} implies that $\operatorname{deg}(\mathcal{C}_{i})  \ge 0 $ for $1\le i \le  N$. By Lemma \ref{lem11},  we have $ \operatorname{deg}(\mathcal{C}^0) \ge 0$. Let $M$ be the total
multiplicity of $\alpha$. Then $\mathcal{C} \cdot \Sigma = M+ \operatorname{deg}(\mathcal{C}_{\diamondsuit}) e + ke = Z_{\alpha} \cdot [\Sigma] = M$. Thus, we have
$k = \operatorname{deg}(\mathcal{C}_{\diamondsuit})= \operatorname{deg}(\mathcal{C}^i)= 0.$ In other words,  $\mathcal{C}$ has no closed components, and  each level of $\mathcal{C}$ has degree $0$.

 To prove the claim,  the argument here is the same as Proposition 4.4 in \cite{NW}.  Let
$S$ be a representative of $[\mathcal{C}_{\diamondsuit}] = Z_{\alpha} + \operatorname{deg}(\mathcal{C}_{\diamondsuit})[\Sigma]$. We have the following energies:
\begin{equation} \label{eq12}
\begin{split}
&\int_{S \cap DE} \Omega_{\varepsilon}  +   \int_{S \cap \mathbb{R}_{s \ge 0} \times Y} d\lambda_{\varepsilon}  = M +  \varepsilon \pi^*H(\alpha)+ \operatorname{deg}(\mathcal{C}_{\diamondsuit}) |e|\\
&\int_{\mathcal{C}^0 \cap DE} \Omega_{\varepsilon}  +   \int_{\mathcal{C}^0 \cap \mathbb{R}_{s \ge 0} \times Y} d\lambda_{\varepsilon}  =  M_{\alpha_0} + \varepsilon \pi^* H(\alpha_0) + \operatorname{deg}(\mathcal{C}^0) |e|\\
& \int_{\mathcal{C}^i \cap \mathbb{R}_{s} \times Y} d\lambda_{\varepsilon}  =2(M_{\alpha_i} - M_{\alpha_{i-1}})  + \varepsilon \left( \pi^*H(\alpha_i) -  \pi^*H(\alpha_{i-1}) \right),
\end{split}
\end{equation}
where $M_{\alpha_i}$ is the total multiplicity of $\alpha_i$,  $\pi^*H(\alpha_i)$ is short for $m_+^i \pi^*H(e_+) +m_-^i \pi^*H( e_-) + \sum^{2g}_{j=1} m_j^i \pi^*H(h_i)$ and $m_{\pm}^i, $ $m_j^i$ are  multiplicities of $e_{\pm}, h_j $ in $\alpha_i$.  Since the energy only
depends on the relative homology class, we have
\begin{equation*}
\begin{split}
\int_{S \cap DE} \Omega_{\varepsilon}  +   \int_{S\cap \mathbb{R}_{s \ge 0} \times Y} d\lambda_{\varepsilon}
=\int_{\mathcal{C}^0 \cap DE} \Omega_{\varepsilon}  +   \int_{\mathcal{C}^0 \cap \mathbb{R}_{s \ge 0} \times Y} d\lambda_{\varepsilon}
 + \sum_{i=1}^N\int_{\mathcal{C}^i \cap \mathbb{R}_{s} \times Y} d\lambda_{\varepsilon}.
\end{split}
\end{equation*}
By Equations (\ref{eq12}), we have
\begin{equation*}
\operatorname{deg}(\mathcal{C}_n) |e|=\operatorname{deg}(\mathcal{C}^0)|e| + \sum_{i=1}^N (M_{\alpha_i} - M_{\alpha_{i-1}}) =\sum_{i=0}^N \operatorname{deg}(\mathcal{C}^i)|e|.
\end{equation*}
The second equality follows from  $\operatorname{deg}(\mathcal{C}^i)|e| =M_{\alpha_i} - M_{\alpha_{i-1}}$ (Proposition 4.4 of \cite{NW}).

 %Recall that $ \operatorname{deg}(\mathcal{C}_n) =0$ because its relative homology class is  $Z_{e_-^M}$. Proposition 4.4 of \cite{NW} also implies that $\operatorname{deg}(\mathcal{C}^i) \ge 0$ for $i=1...N$.  By Lemma \ref{lem11}, we have $\operatorname{deg}(\mathcal{C}^i)  =  0$ for $i=0...N$.  By Lemma \ref{lem6}, $\mathcal{C}^0$ has no closed components.

 %By the ECH index formula and $M =A\cdot B +  d|e|$, we have
%\begin{equation} \label{eq7}
%\begin{split}
%I(\mathcal{C})&= M+m_+-m_- + 2dM +d^2 e +de + d\chi(B)\\
%&\ge  dA\cdot B + d(M-|e|+ \chi(B))\\
%& =  d(2A\cdot B + d|e|-|e|+ \chi(B)) \ge d \ge k.
%\end{split}
%\end{equation}
%Therefore, $k=0$.
\end{proof}

\begin{remark} \label{remark3}
The argument in Lemma \ref{lem6} still works for a general generic cobordism
admissible almost complex structure, because any closed $J$-holomorphic current has
homology class $k[\Sigma]$ for $k \ge 1$. Also, the proof of Lemma \ref{lem11} does not use the conditions \ref{J1} and \ref{J2}.
\end{remark}

\begin{lemma}
The moduli space $\mathcal{M}^J(e_-^M, Z_{e_-^M})$ is a finite set.
\end{lemma}
\begin{proof}
Let $\mathcal{C}_{\infty} =\{\mathcal{C}^{0},...,\mathcal{C}^{N}\}$ be a broken $J$-holomorphic curve which  which is a limit from
a sequence of holomorphic curves  $ \{\mathcal{C}_{n}\}_{n=1}^{\infty}$ in  $\mathcal{M}^J(e_-^M, Z_{e_-^M})$. By Lemma \ref{lem6},   $\mathcal{C}^0$ has no closed components.

Then the rest of the proof is the same as Proposition 3.13 in \cite{CG}. Here we just
sketch the proof. Because the $J$-holomorphic currents in $E$ has nonnegative ECH index,
$I(\mathcal{C}_{\infty}) = I(\mathcal{C}_n) = 0$ implies that $I(\mathcal{C}^i) = 0$ for $0  \le i \le  N$. As a result, $\mathcal{C}^i$ are branched
covers of trivial cylinders. The branched
covers of trivial cylinders have nonnegative Fredholm index (see Exercise 3.14 of \cite{H4}). By Lemma \ref{lem7}, $\mathcal{C}^0$ also has nonnegative Fredholm index.   Since Fredholm index is additive, we  have $\operatorname{ind}(\mathcal{C}^i) = 0$ for
$0\le i\le N$. These branched covers of trivial cylinders can be ruled out by the ECH
partition conditions (see Exercise 3.14 of \cite{H4})). Therefore,  $\mathcal{M}^J(e_-^M, Z_{e_-^M})$ is compact.
The transversality is guaranteed by the automatic transversality theorem \cite{Wen}.

\end{proof}

\begin{comment}
\begin{theorem}
  then for a generic almost complex structure $J$, we have a well--defined homomorphism
\begin{equation*}
ECH^L(E, F_{\varepsilon}^*\Omega_E)_J: ECH^L(Y, \lambda_{\varepsilon}) \to \mathbb{Z}.
\end{equation*}
\end{theorem}
\begin{proof}
The proof is essential the same as before. Let $\alpha $ be an ECH generator. The moduli space $\mathcal{M}_0^J(\alpha)$ is a finite set with orientation can be proof by the same argument in Gerig's paper. Hence, we can define a map in chain level by
\begin{equation*}
ECC^L(E, F_{\varepsilon}^*\Omega_E)_J\alpha=\#\mathcal{M}^J_0(\alpha).
\end{equation*}
To see this is  a chain map, let $\mathcal{C}$ be a  broken holomorphic curves come from the limit of curves in $\mathcal{M}^J_1(\alpha)$. Then $\mathcal{C}$ is either a curve in  $\mathcal{M}^J_1(\alpha)$ or consists of
\begin{itemize}
\item
A holomorphic curve with $I=1$ in the top level;
\item
A holomorphic curve with $I=0$ in the cobordism level;
\item
Connectors in the middle level.
\end{itemize}
We can apply HT's gluing argument, $ECC^L(E, F_{\varepsilon}^*\Omega_E)_J$ is a chain map.
\end{proof}
\end{comment}
%\begin{lemma}
%If $I(\mathcal{C}) =1$, then $\mathcal{C} =\mathcal{C}_{emb} \cup \mathcal{C}_{spec}$, where $\mathcal{C}_{emb}$ is embedded holomorphic current with $I=ind=1$ and $\mathcal{C}_{spec}$ consists of special holomorphic planes.
%\end{lemma}
\paragraph{Uniqueness}
We show that  $MC_{e_-}$ is the unique $J$-holomorphic current in the moduli space $\mathcal{M}^J(e_-^M, Z_{e_-^M})$. The energy constraint argument in \cite{GHC1} does not work in the current situation. We use the argument as in Lemma \ref{lem7} instead.

To this end, we need to apply R. Siefring's intersection theory for punctured $J$-holomorphic
curves.  In  \cite{RS},  Siefring defines a intersection pairing   $C\bullet C'$ for punctured holomorphic curves, where $C$ and $C'$ are simple holomorphic curves.
Here we do not use the precise definition of $\bullet $, we only need  to know the following facts:
\begin{enumerate} [label=\textbf{ F.\arabic*}]
  \item \label{F1}
 The  intersection pairing  $C\bullet C'$ is invariant under homotopy as cylindrical asymptotic
maps.
 \item\label{F2}
 (Theorem 2.2 of \cite{RS})
 If   $C$ and $C'$ are  distinct, then $C\bullet C' \ge 0$.
  \item \label{F3}
   (Theorem 2.3 of \cite{RS})
 In the case that $C=C'$, we have the following adjunction
formula:
 \begin{equation*}
C\bullet C=2(\delta(C)+ \delta_{\infty}(C)) +\frac{1}{2} (2g(C)-2 + \operatorname{ind}(C) + \# P_{even})+ (\bar{\sigma}(C)-\#P),
\end{equation*}
where $P$ denote the set of punctures, $P_{even}$ is the set of punctures which are asymptotic to  Reeb orbits with even Conley-Zehnder index, $\delta_{\infty}(C)$ is an algebraic count of ``hidden'' singularities at the infinity, and $\bar{\sigma}(C)$ is the spectral covering
number. We refer readers to \cite{RS} for the precise definition of  $\delta_{\infty}(C)$ and $\bar{\sigma}(C)$.
According to the definition, if all the ends of $C$ are asymptotic to distinct simple
orbits, then $\delta_{\infty}(C)$ and $\bar{\sigma}(C)-\#P$ vanish.
\end{enumerate}

\begin{lemma} \label{lem3}
%Fix a relative class $A \in H_2(E, Y, \mathbb{Z})$ represented by $Z_{e_-^M}$.
%Then the  moduli space $\mathcal{M}^J_0(\alpha, Z)$ is empty if $\alpha \ne  e_-^M$. If $\alpha = e_-^M$,
 The moduli space    $\mathcal{M}^J(e_-^M, Z_{e_-^M})$ only consists of one element.
\end{lemma}
\begin{proof}
Note that $MC_{e_-} \in   \mathcal{M}^J(e_-^M, Z_{e_-^M})$.  Moreover,    $I(MC_{e_-}) = I(Z_{e_-^M})=0$. Thus, the  moduli space is nonempty.

 Let $\mathcal{C}=\sum_a d_a C_a \in \mathcal{M}^J(e_-^M, Z_{e_-^M}) $.  By Lemma \ref{lem6}, $\mathcal{C}$ has no closed components, and  $\operatorname{deg}(C_a)=0$ for any $a$. %Lemma \ref{lem1} implies that $I(\mathcal{C} ) =0$ if and only if $\alpha=e_-^M$.
 By  (\ref{eq5}) and (\ref{eq20}), we have $I(C_a)=\operatorname{ind}( C_a)=\delta(C_a)=0$ for every $a$. Lemma \ref{lem2} forces $M_a=1$, $g(C_a)=0$ and $h(C_a) =e_+(C_a)=0$. In sum,   $C_a$ are special $J$-holomorphic planes.

By our choice of $J$, the fiber $C_{e_-} $ is a holomorphic plane with $I(C_{e_-}) =\operatorname{ind}(C_{e_-})
=0$.  By \ref{F3}, we have $C_{e_-} \bullet C_{e_-}=-1$. If there exists another special $J$-holomorphic plane $C_a $ other than $C_{e_-}$, then $C_{e_-} \bullet C_a \ge 0$ by \ref{F2}. Note that $C_a$ is homotopy  to $C_{e_-}$ as a  asymptotic cylindrical map because $\pi_E(C_{e_-} -C_a)$ is trivial in $\pi_2(\Sigma)$. Therefore, \ref{F1} implies that $0 \le C_a \bullet C_{e_-}=C_{e_-}\bullet C_{e_-}=-1$. We get a contradiction.
 \end{proof}

\paragraph{$(L, \delta)$-flat approximation}
Fix an admissible almost complex structure $J \in \mathcal{J}(Y, \lambda_{\varepsilon})$.  To ensure that the ECH generators are in a one-to-one correspondence with the solutions
to the Seiberg-Witten equations,  we need to perturb $(\lambda_{\varepsilon}, J)$ such that it has certain  standard forms ((2-11) of \cite{T2}) in  $\delta$-neighborhoods of Reeb orbits with action less than $L$. Moreover,  the modifications do not change  the ECH chain complex.    The result of the perturbation   of $(\lambda_{\varepsilon}, J)$   is called a \textbf{$(L, \delta)$-flat approximation}, and it was introduced by Taubes \cite{T2}.  The $(L, \delta)$-flat approximation closes to the original one in $C^0$ topology.

Suppose that $(\lambda, J)$ is $(L, \delta)$-flat, then we have a canonical isomorphism (Theorem 4.2 of \cite{T2})
\begin{equation}
  \Psi: ECC_*^L(Y, \lambda, \Gamma) \to CM^{-*}_L(Y, \lambda, \mathfrak{s}_{\Gamma})
\end{equation}
 between the ECH chain complex and the Seiberg-Witten chain complex, where $CM^{-*}_L(Y, \lambda, \mathfrak{s}_{\Gamma})$
is the Seiberg-Witten chain complex defined in Section 2.3 of \cite{HT}, and $\mathfrak{s}_{\Gamma}$ is the $Spin^c$ structure such
that $c_1(\mathfrak{s}_{\Gamma}) = c_1(\xi) + 2\operatorname{PD}(\Gamma)$.

 %By Lemma 3.6 of \cite{HT}, there exists a preferred homotopy $\{(\lambda_{\varepsilon}^t, J^t)\}_{t\in [0, 1]}$ such that  $(\lambda_{\varepsilon}^0, J^0) = (\lambda_{\varepsilon}, J)$ and $(\lambda_{\varepsilon}^1, J^1)$ is  $(L, \delta)$-flat. Moreover,  $(\lambda_{\varepsilon}^t, J^t)$ is independent of $t$ outside the  $\delta$-neighborhoods of Reeb orbits with action less than $L$.

 To apply the correspondence between solutions to Seiberg-Witten equations and $J$-holomorphic
curves, we need to modify $(\Omega_{\varepsilon}, J)$ such that the induced contact form and
almost complex structure are ($L$-$\delta$)-flat over the cylindrical end. Fix  $\epsilon>0$   and $\delta>0$. By the same construction as
in Section 6.3 of  \cite{HT}, we can define a new symplectic form $\Omega_{\delta} $ and  an almost complex structure $J_{\delta}$ such that
\begin{enumerate}[label=\textbf{ C.\arabic*}]
  \item\label{C1}
  $\Omega_{\delta} =d(e^s \lambda_{\delta})$ on the cylindrical end $[0,\infty) \times Y$. $J_{\delta}$ is a cobordism admissible almost complex structure with respect to    $\Omega_{\delta}$.  Furthermore, $(\lambda_{\delta}, J_{\delta})$ is a ($L$-$\delta$)-flat approximation of $(\lambda_{\varepsilon}, J)$.
   \item \label{C2}
    $(\Omega_{\delta}, J_{\delta}) =(\Omega_{\varepsilon}, J)  $ on the region $s \le -\epsilon$.
  \item \label{C3}
  $J_{\delta}$ satisfies estimates $|J_{\delta}-J|_{C^0} \le c_0\delta$ and $|J_{\delta}-J|_{C^1} \le c_0$ for some $\delta$-independent
constant $c_0$.
\end{enumerate}

The following lemma tell us that the conclusion in Lemma \ref{lem3} is still true if we
replace         $(\Omega_{\varepsilon}, J) $ by $(\Omega_{\delta}, J_{\delta}) $.

\begin{lemma}\label{lem10}
For sufficiently small $\delta>0$, $\mathcal{M}^{J_{\delta}} (e_-^M, Z_{e_-^M}
) $ only consists of one element.
\end{lemma}
\begin{proof}
From the proof of Lemma \ref{lem3}, we know that $\mathcal{M}^{J_{\delta}} (e_-^M, Z_{e_-^M}
) $  only consists of unbranched
cover of special $J_{\delta}$-holomorphic plane. Moreover, the special $J_{\delta}$-holomorphic
plane is unique provided that it exists.

It suffices to show that the special $J_{\delta}$-holomorphic plane exists for small $\delta$. By
(2-11) of \cite{T2}, $J_{\delta} = J$ along $\mathbb{R}_{\ge 0} \times e_-$.  Since $J_{\delta} = J$ for $s \le -\epsilon$, $C_{e_-}$ is $J_{\delta}$-holomorphic
except the part in the small region $[-\epsilon, 0] \times U$, where $U$ is a $\delta$-neighborhood of $e_-$. $C_{e_-}$ is nearly $J_{\delta}$-holomorphic in the sense that it satisfies (A-15) in \cite{T2}. Therefore, we
can apply the same argument in Appendix A of \cite{T2} to deform $C_{e_-}$ so that it becomes $J_{\delta}$-holomorphic. We sketch the main argument as follows.

Let $N_C $ denote the normal bundle of $C_{e_-}$ and $e_C : N_C \to E$ denote an exponential
map. Let $\eta$ be a section of $N_C$ which is asymptotic to zero as $s \to \infty$. Then $e_C(\eta)$ is $J_{\delta}$-holomorphic if and only if $\eta$ satisfies the following equation:
\begin{equation}\label{eq21}
  D_C\eta + \mathfrak{p}_1 \eta +  \mathfrak{p}_2\nabla_C \eta  +  \mathfrak{R}_0 \eta +  \mathfrak{R}_1 \nabla_C\eta + \mathfrak{p}_0=0.
\end{equation}
Here $\mathfrak{p}_0 = \bar{\partial}_{J_{\delta}} C_{e_-}$ satisfying $|\mathfrak{p}_0|_{L^2} \le c_0 \delta$ due to the discussion in the last paragraph and
the estimates of $J_{\delta}$. The operators $\mathfrak{p}_1$, $\mathfrak{p}_2$ come from the linearization of $(J_{\delta}-J) \circ de_C(\eta)\circ
j$ and they satisfy estimates $|\mathfrak{p}_1| \le c_0|J -J_{\delta}|_{C^1} \le  c_0 $ and $|\mathfrak{p}_2| \le  c_0|J -J_{\delta}|_{C^0} \le  c_0 \delta $.  $\mathfrak{R}_0$
and $\mathfrak{R}_1$ are high order terms and they satisfy (A-17) in \cite{T2}.

If the operator $D_C +  \mathfrak{p}_1   +  \mathfrak{p}_2\nabla_C : L_1^2(N_C ) \to L^2(N_C \otimes T^{0, 1}C) $ is invertible, then we can apply the contraction mapping argument in Lemma A.3 of \cite{T2} to find a solution $\eta$
satisfying (\ref{eq21}). Then, $\mathcal{M}^{J_{\delta}} (e_-^M, Z_{e_-^M})  $ is nonempty and it contains exactly one element.

Since $|\mathfrak{p}_2| \le c_0 \delta$, it suffices to prove that  $D_C +  \mathfrak{p}_1 $ is invertible. Note that $D_C +  \mathfrak{p}_1 $ is
a CR-operator which is asymptotic to the operator $\frac{i}{2} \frac{d}{dt} + \frac{\theta}{2}$ with odd Conley-Zehnder.
Because the  ($L$-$\delta$)-approximation does not change the rotation number of the Reeb orbit
(see (2-11) of \cite{T2}), we still have $\operatorname{ind}(D_C + \mathfrak{p}_1) = 0$. By the automatic transversality
theorem \cite{Wen},  $ D_C + \mathfrak{p}_1 $  is invertible.
\end{proof}

\begin{proof} [Proof of Theorem \ref{thm1}]
Let $A \in H_2(DE, Y, \mathbb{Z})$ be the relative class represented by $[Z_{e_-^M}].$ Recall that $Z_{e_-^M}$ is the  unique  relative homology class $Z$ in $H_2(DE, e_-^M)$ such that $[Z]=A$.

Let $(\Omega_{\delta}, J_{\delta}) $ be the  pair satisfying \ref{C1}, \ref{C2} and \ref{C3}.  Then we have  a bijection (\ref{eq16}) between the ECH generators  and the gauge  classes  of   the Seiberg-Witten solutions.   Let $\mathfrak{c}_{e_-^M} =\Psi(e_-^{M})$. The ECH cobordism map is defined by (Definition 5.9 of \cite{HT})
\begin{equation*}
ECC^L(DE, \Omega_{\varepsilon}, A)({e_-^M}) =\#\mathfrak{M}(\mathfrak{c}_{e_-^M}, \mathfrak{s}_A),
\end{equation*}
 where $ \mathfrak{M}(\mathfrak{c}_{e_-^M}, \mathfrak{s}_A)$ is the moduli space of solutions to  the Seiberg-Witten equations  on $E$  which are asymptotic to   $\mathfrak{c}_{e_-^M}$ (see (4.15) of \cite{HT}), $\mathfrak{s}_A$ is the  $Spin^c$ structure such that $c_1(\mathfrak{s}_A) =c_1(K_{DE}^{-1}) +2\operatorname{PD}_{DE}(A)$ and $K_{DE}^{-1}$ is the canonical line bundle. By Theorem 4.2 of \cite{CG2} and  Lemma \ref{lem3},  we have $$\#\mathfrak{M}(\mathfrak{c}_{e_-^M}, \mathfrak{s}_A) = \# \mathcal{M}^{J_{\delta}}(e_-^M, Z_{e_-^M}) =1.$$
 By Lemma 4.6 and Corollary 6.2 of \cite{NW}, the differential vanishes on $ECC^{L_{\varepsilon}}(Y, \lambda_{\varepsilon})$.
Then  $e_-^M$ is a cycle,  and we have $ECH^L(DE, \Omega_{\varepsilon}, A)([{e_-^M}]) =1.$

Let $\beta =e_+^{n_+} h_1^{n_1},..., h_{2g}^{n_{2g}} e_-^{n_-}$ be an ECH generator with $N = M-d|e|$ for $d \ge 1$, where $N=n_+ \sum_{i=1}^{2g}n_i + n_-$.
To see $ECH^L(DE, \Omega_{\varepsilon}, A)( [\beta]) =0,$ by the holomorphic curve axioms (see Theorem 1.9 of \cite{HT}), it suffices to show that the moduli space $\overline{\mathcal{M}^J(\beta, Z_{\beta, A})}$ is empty, where $Z_{\beta, A} \in H_2(DE, \beta)$  is the unique relative homology class determined by $A$. Let $\mathcal{C}= \mathcal{C}_{\diamondsuit} + k\Sigma$ be a holomorphic current in this  moduli space, where $ \mathcal{C}_{\diamondsuit}$ has no closed component and $k \ge 0$. Then
\begin{equation*}
  \mathcal{C} \cdot \Sigma = M-|e| + \operatorname{deg}( \mathcal{C}_{\diamondsuit})e +ke =A \cdot \Sigma =Z_{e_-^M} \cdot \Sigma =M.
\end{equation*}
Then $\operatorname{deg}( \mathcal{C}_{\diamondsuit}) +k=-d\le -1$. This  contradicts with Lemma \ref{lem11} (also see Remark \ref{remark3}).
\end{proof}

\subsection{Proof of Theorem \ref{thm1}}
To prove Theorem \ref{thm1} and Theorem \ref{thm4}, we first need to compute the $U$ map for some
ECH generators. The computations are parallel to Lemma 4.6 of \cite{NW}.

Let $\mathbf{z} = \{z_1, ..., z_k\}$ be $k$ distinct marked points in $Y$ away from the Reeb orbits.
Let $\mathcal{M}_i^{J}(\alpha, \beta)_{\mathbf{z}}$  denote the moduli space of ECH index-$i$ $J$-holomorphic currents passing
through the marked points $\{0\}\times \mathbf{z} \in \mathbb{R} \times Y$. By the same argument as in Lemma 2.6 of
\cite{HT3},  $\mathcal{M}_{2k}^{J}(\alpha, \beta)_{\mathbf{z}}$ is a finite set for a generic almost complex structure. In the case that
$k = 1$, the counting of this moduli space is used to define the $U$ map. By the similar
argument in Proposition 3.25 of \cite{CG}, we can define the $U^k$  at the chain level by
\begin{equation} \label{eq22}
 U^k \alpha : =\sum_{\beta}\# \mathcal{M}_{2k}^{J}(\alpha, \beta)_{\mathbf{z}} \beta
\end{equation}

\begin{lemma} \label{lem12}
Fix a positive integer $k$, and let $J$ be a generic admissible almost complex structure. Then the moduli space $\mathcal{M}_{2k}^J(e_+^k, e_-^k)_{\mathbf{z}}$ consists of exactly one element, which is made up of $k$ distinct index-$2$ holomorphic cylinders passing through the marked point $\mathbf{z}$. As a result, we have
 $\langle U^k e_+^k, e_-^k \rangle=1$.
\end{lemma}
\begin{proof}
 Let $\mathcal{C} =\sum_a d_a C_a$ be a holomorphic current in  $\mathcal{M}_{2k}^J(e_+^k, e_-^k)_{\mathbf{z}}$. Note that $I(\mathcal{C})=I(e_+^k, e_-^k) =2k$ by by the index formula (\ref{eq8}).

 Since $J$ is generic,     if $C_a$ passes through  $l_a \ge 0 $ marked points, then $\operatorname{ind}(C_a) \ge 2l_a$. By Hutchings's ECH inequality (see  Theorem 4.15 and Theorem 5.1 of  \cite{H2} for example), we have
  \begin{equation*}
  \begin{split}
    2k=I(\mathcal{C}) &\ge \sum_a d_aI(C_a) + 2\sum_{a \ne a '} d_a d_{a'}C_a \star C_{a'}\\
     &\ge \sum_a d_a \left( \operatorname{ind}(C_a) + 2 \delta(C_a)\right) \\
     &\ge \sum_{a, l_a \ge 1} 2l_a+  \sum_{a, l_a \ge 1}  2(d_a-1)  l_a    +  \sum_{a }   2 d_a    \delta(C_a) .
  \end{split}
  \end{equation*}
Since $\sum_a l_a=k$, we must have $d_a=1$,  $I(C_a) = \operatorname{ind}(C_a) =2l_a$ provided that $l_a \ge 1$ and $C_a$ is a trivial cylinder provided that  $l_a=0$. Furthermore,  $C_a \star C_{a'} = 0$, i.e., ${C_a}$ are pairwise disjoint.  By  Proposition 4.4 of \cite{NW}, the  degree of $ {C_a}$ are nonnegative and additivity. Hence, we must have $\operatorname{deg}(C_a) =0$ because their sum is  $\operatorname{deg}(\mathcal{C}) = (k-k)/|e|=0$. Consequently, we have $C_a \in \mathcal{M}^J(e_+^{m_a}, e_-^{m_a})$.

By Nelson and Weiler's index formula \ref{eq8} and the ECH partition condition, we have
  \begin{equation}\label{eq24}
  \begin{split}
    &2l_a= \operatorname{ind}(C_a) = 2g(C_a) -2 + 4m_a, \\
   & 2l_a= I(C_a) =  2m_a.
  \end{split}
  \end{equation}
Hence, we have $g(C_a)=0$ and $l_a =m_a =1$, i.e., $C_a$ is a holomorphic cylinder from $e_+$ to $e_-$ passing a marked point. Consequently, there are $k$ holomorphic cylinders and there are  no trivial cylinders.  By Proposition 4.7 of \cite{NW}, there is a bijection between $\mathcal{M}_{2k}^J(e_+^k, e_-^k)_{\mathbf{z}}$ and the moduli space of Morse flow lines passing through the marked points. For each marked point, there is exactly one index $2$ Morse flow line passing through it. Therefore, $\mathcal{M}_{2k}^J(e_+^k, e_-^k)_{\mathbf{z}}$ contains exactly one element and  $\langle U^k e_+^k, e_-^k \rangle= \#\mathcal{M}_{2k}^J(e_+^k, e_-^k)_{\mathbf{z}}=1$.
\end{proof}

We have the following parallel version of Lemma \ref{lem12}  for $J$-holomorphic currents
in $E$ with points constraint. Let $\mathcal{M}^J (\alpha,Z)_{\mathbf{z}}$ denote the moduli space of $J$-holomorphic
currents in $E$ passing through the marked points $\mathbf{z}$.

\begin{lemma} \label{lem13}
Let $J$ be a generic cobordism admissible almost complex structure (not
necessarily satisfying the conditions \ref{J1} and \ref{J2}) and $\mathbf{z} = \{z_1, ..., z_k\} \in  \operatorname{Int}(DE)$. Then
there is a k disjoint union of embedded $J$-holomorphic planes in  $\mathcal{M}^J (e^k_+,Z_{e^k_+})_{\mathbf{z}}$.
\end{lemma}
\begin{proof}
Let $\beta=e_+^{n_+} h_1^{n_1},...,h_{2g}^{n_2g} e_-^{n_-}$.  Suppose that $\langle U^ke^k_+, \beta \rangle= 1$. By definition of the U map,   there exists  a $J$-holomorphic
current $\mathcal{C} \in \mathcal{M}_{2k}^J(e_+^{k}, \beta)_{\mathbf{z}}$. By Proposition 4.4 of \cite{NW}, the degree of $\mathcal{C}$ is
$\operatorname{deg}(\mathcal{C}) =\frac{k-N}{|e|} \ge 0$, where $N = n_+ +
\sum_{i=1}^{2g}n_i + n_-.$ Then $N = k - \operatorname{deg}(\mathcal{C}) |e| \le k$. If
$N = k$, then $\beta =e_-^k$  by the index formula (\ref{eq8}).

Let $A = [Z_{e_-^k}] \in H_2(DE, Y, \mathbb{Z})$. By Lemma \ref{lem12}, Theorem \ref{thm2} and the above discussion, we have
\begin{equation*}
  ECH^L(DE, \Omega_{\varepsilon}, A)(U^k[e^k_+]) =  ECH^L(DE, \Omega_{\varepsilon}, A)([e^k_-])=1.
\end{equation*}
Note that $Z_{e_+^k, A} =Z_{e_+^k}$. By the holomorphic curve axiom \cite{HT}, we have a broken $J$-holomorphic
current $\mathcal{C} =\{\mathcal{C}^0, ..., \mathcal{C}^N\} \in \overline{\mathcal{M}^J(e_+^k, Z_{e_+^k})_{\mathbf{z}}}$, where $\mathcal{C}^0 \in \mathcal{M}^J(\alpha_0)_{\mathbf{z}}$, $\mathcal{C}^i \in\mathcal{M}^J(\alpha_i, \alpha_{i-1})$ ($1\le i \le N$), and $\alpha_N=e_+^k$.  By Lemma \ref{lem6} and Remark \ref{remark3}, there is no closed components and the
degree of each level is zero.

Write $\mathcal{C}^0 =\sum_a d_a C_a$. By Lemma \ref{lem1}, $I(d_aC_a) \ge 0 $ because its degree is zero. Since
$J$ is generic, $\operatorname{ind}(C_a) \ge  2l_a$ provided that $C_a$ passes through $l_a\ge 0$ marked points. By definition,
$C_a \star C_a  \ge  0$  if $C_a$ passes at least one marked point. By (\ref{eq11}) and (\ref{eq20}),
we have
\begin{equation}\label{eq23}
   \begin{split}
   2k=I(\mathcal{C}) =\sum_{i=0}^NI(\mathcal{C}^i) &\ge I(\mathcal{C}^0)   \\
   & \ge \sum_{a, l_a \ge 1} d_a I(C_a) +  \sum_{a, l_a =0} I(d_aC_a)  +2\sum_{a\ne a'} C_a \star C_{a'}\\
   &\ge 2k + \sum_{a, l_a \ge 1} 2(d_a-1)l_a +2\sum_{a, l_a \ge 1} \delta(C_a).
   \end{split}
\end{equation}
As a result, if $l_a \ge 1$, then  $d_a=1$. Moreover, $I(\mathcal{C}^0) =2k$ and $I(\mathcal{C}^i) =0$ for $i \ge 1$.  By  Proposition 3.7 of \cite{H4}, we have  $\mathcal{C}^i=k (\mathbb{R} \times e_+)$ for $i \ge 1$ and $\alpha_0 =e_+^k$.  The inequality (\ref{eq23}) also implies that $\{C_a\}$ are
pairwise disjoint,  and $C_a$ is embedding provided that  $l_a \ge  1$.

Let $d_aC_a$ be a component in the summation
$\sum_{a, l_a=0} I(d_aC_a)$. The degree of $C_a$
is zero and the ends of $C_a$ are asymptotic to covers of $e_+$. By Lemma \ref{lem1}, we have
$I(d_aC_a) \ge 2$. This contradicts (\ref{eq23}). Therefore, each component of $\mathcal{C}^0$ passes at least
one marked point. For each $C_a$, by Lemma \ref{lem1} and Lemma \ref{lem2}, the ECH index and
Fredholm index of $C_a$ satisfy (\ref{eq24}). Then we obtain $l_a = m_a = 1$ and $g(C_a) = 0$.
Thus, $C_a$ is a $J$-holomorphic plane and it passes exactly one marked point.

 In sum, $\mathcal{C}^0 = \sum_{a=1}^k C_a$  is a $k$ disjoint union of embedded $J$-holomorphic planes in ${\mathcal{M}^J(e_+^k, Z_{e_+^k})_{\mathbf{z}}}$.
 \end{proof}

 \begin{remark}
 By the correspondence between holomorphic curves and the solutions to
the Seiberg-Witten equations (see  Theorem
4.2 of   \cite{CG2}), we should have $\#{\mathcal{M}^J(e_+^k, Z_{e_+^k})_{\mathbf{z}}}=1$. But  for our purpose,  we
do not  need this stronger result.
 \end{remark}

 \begin{proof}[Proof of Theorem \ref{thm1}]
 Suppose that we have a symplectic embedding $$\varphi : (\sqcup_{i=1}^k (B^4(r_i), \omega_{\operatorname{std}}) \hookrightarrow
(DE, \Omega).$$ Let  $\mathbf{z}= \{z_1, ..., z_k\} \in DE$ be the image of the center of the balls. Since $\Omega_{\varepsilon} =\Omega$ except in a small collar neighbourhood of $Y$.  For sufficiently small $\varepsilon>0$, we can choose
such a collar neighbourhood such that it is disjoint from the image of $\varphi$. Thus, we regard
$\varphi$ as a symplectic embedding of  $(\sqcup_{i=1}^k (B(r_i), \omega_{\operatorname{std}})$ into $(DE,\Omega_{\varepsilon})$.

 Let $J_0$ be the standard almost complex structure on $\sqcup_{i=1}^k B^4(r_i)$. Extend $\varphi_*(J_0)$ to be a
generic cobordism admissible complex structure $J$. By Lemma \ref{lem13}, there is a $k$ disjoint  union of
embedded $J$-holomorphic planes $\mathcal{C}$ passing through $\mathbf{z}$. By the standard monotonicity
lemma (see Theorem 5.2.1 of \cite{Wen2}) and (\ref{eq12}), we obtain
\begin{equation*}
  \sum_{i=1}^k r_i \le \int_{\mathcal{C} \cap \varphi(\sqcup_{i=1}^k (B^4(r_i))} \Omega_{\varepsilon} \le \int_{\mathcal{C} \cap DE} \Omega_{\varepsilon} + \int_{\mathcal{C} \cap (\mathbb{R}_{s\ge0} \times Y)} d\lambda_{\varepsilon} \le k+O(\varepsilon).
\end{equation*}
Consequently, we get the conclusion of Theorem \ref{thm1} by taking $\varepsilon \to 0$.
 \end{proof}

\section{Proof of Theorem \ref{thm3} and Theorem \ref{thm4}}
To prove Theorem \ref{thm3} and Theorem \ref{thm4}, we need the following spectrality of the ECH
spectrum.

\begin{lemma}\label{lem14}
Suppose that $(Y, \lambda)$ is $L_0$-nondegenerate and $c_k(Y, \lambda) < L_0$. Then there
exists a null homologous ECH generator  $\alpha$ such that $c_k(Y, \lambda) = \mathcal{A}_{\lambda}(\alpha)$ and $\operatorname{gr}(\alpha) = 2k$.
\end{lemma}
\begin{proof}
Let $\sigma \in ECH_{2k}(Y, \lambda, 0)$ with $U^k\sigma =[\emptyset]$. Define \begin{equation*}
  c_{\sigma}(Y, \lambda):=\inf\{L \vert \sigma \in \mbox{image of } i_L: ECH^L_{2k}(Y, \lambda, 0) \to ECH_{2k}(Y, \lambda, 0) \}.
\end{equation*}
By definition, we have $c_{k}(Y, \lambda)=\min\{c_{\sigma}(Y, \lambda) \vert \sigma \in ECH_{2k}(Y, \lambda, 0), U^k \sigma=[\emptyset]\}$.

Assume that $c_{\sigma}(Y, \lambda)<L_0$. By definition, we can rewrite it as
\begin{equation}\label{eq17}
  c_{\sigma}(Y, \lambda):=\inf\{\mathcal{A}_{\lambda}(x) \vert  x=\sum a_i \alpha_i  \in  ECC^{L_0}_{2k}(Y, \lambda, 0), \partial x= 0,  \mbox{ and } i_{L_0}([x]_{L_0}) =\sigma  \},
\end{equation}
where  $\mathcal{A}_{\lambda}(x) :=\max\{\mathcal{A}_{\lambda}(\alpha_i) \vert a_i \ne 0\}$ and  $[x]_{L_0}$ denotes the class in  $ECH^{L_0}_{2k}(Y, \lambda, 0)$ represented by $x$. Since $\lambda$ is $L_0$-nondegenerate, there are only finitely many ECH
generators with $\mathcal{A}_{\lambda} < L_0$. We can find $\alpha$ such that $c_{\sigma}(Y, \lambda) = \mathcal{A}_{\lambda}(\alpha)$ and $\operatorname{gr}(\alpha) = 2k.$
Take $\sigma$ to the class such that $c_k(Y, \lambda) = c_{\sigma}(Y, \lambda)$; then we obtain conclusion of the
lemma.

\end{proof}
\subsection{Sphere case}
In this subsection,  we assume that  $\Sigma =\mathbb{S}^2$.  It is well known that    the diffeomorphism type of $Y$ is the lens space  $L(|e|, 1)$. The ECH group of $Y$ (as an $\mathbb{F}$ module) has  been computed by  Nelson and  Weiler (Example 1.3) \cite{NW}.  But we still need  to know the $U$ module structure of $ECH(Y,\lambda_{\varepsilon}, 0)$ by using  Taubes's isomorphism ``ECH=HM'' \cite{T2, Te2,Te3,Te4,Te5} and the  computations of P. Kronheimer,  T. Mrowka,  P. Ozsv\'{a}th and Z. Szab\'{o} in \cite{KMOS}.
\begin{prop} \label{lem8}
The ECH of the lens space $Y \cong L(|e|, 1)$ is
 \begin{equation} \label{eq9}
ECH_*(Y, \lambda, 0 ) =
\begin{cases}
\mathbb{F} , & *=2k \mbox{ and } k \ge 0,\\
0, & else,
\end{cases}
\end{equation}
where the $\mathbb{Z}$ grading is defined by (\ref{eq10}).
Moreover, $U: ECH_{2k}(Y, \lambda, 0 )  \to ECH_{2k-2}(Y, \lambda, 0 )$ is an isomorphism for   $k \in \mathbb{Z}_{\ge 1}$.
   Also,  $ECH_0(Y, \lambda, 0) $ is spanned by $[\emptyset]$.
\end{prop}
\begin{proof}
The isomorphism (\ref{eq9}) is just the sphere case of  Theorem 1.1 in \cite{NW}. It remains to show that the $U$ map is an isomorphism.

By Taubes's series papers \cite{T2, Te2,Te3,Te4,Te5}, we have a canonical isomorphism $ECH_*(Y, \lambda, 0) \cong \widehat{HM}^{-*}(Y, \mathfrak{s}_{\xi}) $ as an $U$-module, where  $ \widehat{HM}^{-*}(Y, \mathfrak{s}_{\xi})$ is the hat-version of the Seiberg-Witten Floer cohomology.   Since $L(|e|, 1)$ admits a metric with positive scalar curvature, by Proposition 2.2 and Corollary 2.12 of \cite{KMOS}, we have  the $\mathbb{F}[U]$ module isomorphism  $$\widehat{HM}(Y, \mathfrak{s}_{\xi}) \cong \mathbb{F}[U^{-1}, U]/ \mathbb{F}[U]. $$ Therefore, $U$ is an isomorphism when  the grading is at least two.
\end{proof}

\begin{remark}
We remark that the computation of $ECH(L(|e|, 1), \lambda, \Gamma)$ holds similarly for nonzero
homology class $\Gamma$ (see Corollary 3.4 of  \cite{KMOS}).
\end{remark}

%\begin{lemma}
%%If $ECH^{L_0}(X, \Omega): ECH_{2k}^{L_0}(Y, \lambda) \to \mathbb{F}$ is nonzero, then $c_k(Y, \lambda) \le L_0$.
% Let $d$ be the maximal integer such that $A_{\lambda_{\varepsilon}} (e_-^{d|e|}) <L_{\varepsilon}$. Then for any $k\le gr(e_-^{d|e|})$, the cobordism map
%\begin{equation*}
%ECH(X, \Omega_{\varepsilon}) : ECH^{L_{\varepsilon}}(Y, \lambda_{\varepsilon}) : \to \mathbb{F}
% \end{equation*}
%  is an isomorphism.
%\end{lemma}
%\begin{proof}
%%Let $\sigma \in ECH_{2k}^{L_0}(Y, \lambda)$. Since  $1=ECH^{L_0}(X, \Omega) (\sigma) = ECH(X, \Omega) \circ i_{L_0}(\sigma)$. Therefore, $i_{L_0} (\sigma) \ne 0$.
%We have the following diagram
%$$\begin{CD}
%				ECH_{2k_0}^L(Y, \lambda_{\varepsilon}, 0) @> ECH^L(X,\Omega_{\varepsilon}) >> \mathbb{F}\\
%				@VV U^{k} V @VV  V\\
%				ECH_{2k_0-2k}^L(Y, f\lambda_{\varepsilon}, 0)@> ECH^L(X,\Omega_{\varepsilon})>> \mathbb{F}  \\
%			\end{CD}$$
%\end{proof}

%\begin{lemma}
%$c_{k}(Y, \lambda_{\varepsilon}) =2|e|d $, where $k =d+ \frac{d(d-1)}{2}|e|$.
%\end{lemma}
\begin{proof}[Proof of Theorem \ref{thm3}]
By Proposition \ref{lem8}, we have a sequence $\sigma_k \in ECH_{2k}(Y, \lambda,0)$ such
that $U^k\sigma_k = [\emptyset]$ and $U\sigma_k = \sigma_{k-1}$. By the ``volume property" of ECH spectrum (Theorem
1.3 of \cite{CHR}) which is discovered by D. Cristofaro-Gardiner, M. Hutchings and V.G.B.
Ramos, we know that $c_k(Y, \lambda)$ is finite for each $k$.

  For each $k$, since $c_k(Y, \lambda)$ is finite, there is a constant $L_k > 0$ such that $c_k(Y, \lambda) < L_k$.   Then for $0 < \varepsilon \ll 1$, we still have $c_k(Y, \lambda_{\varepsilon}) < L_k$. Recall that $\lim_{\varepsilon \to \infty} L_{\varepsilon} = \infty$. Thus, $L_{\varepsilon} > L_k$ for sufficiently small $\varepsilon$.   Since $\lambda_{\varepsilon}$ is $L_{\varepsilon}$-nondegenerate, by Lemma \ref{lem14}, we have
an orbit sets $\alpha=e_+^{m_+}e_-^{m_-}$  such that $c_k(Y, \lambda_{\varepsilon}) = \mathcal{A}_{\lambda_{\varepsilon}}(\alpha)$, $\operatorname{gr}(\alpha) = 2k$ and   $m_-+m_+=0 \mod |e|$. The last condition is equivalent to $[\alpha]=0$.

   In \cite{NW}, Nelson and  Weiler show that there is a bijection between  the nonnegative integers $k$ and the pairs $(m_-, m_+)$  satisfying $m_-+m_+=0 \mod |e|$.    Therefore, the unique ECH generator $e_+^{m_+}e_-^{m_-} \in ECC_{2k}^{L_{\varepsilon}}(Y, \lambda_{\varepsilon}, 0)$ is
characterized by the following relation
  \begin{equation} \label{eq13}
  \begin{split}
     & m_-+m_+ =d|e| \mbox{ for some positive integer $d$}, \\
      & 2k= \operatorname{gr}(e_-^{m_-}e_+^{m_+}) =2d+ d^2|e| + m_+ - m_-.
  \end{split}
  \end{equation}
By  (\ref{eq13}), we get
\begin{equation}\label{eq25}
2d+d|e|(d-1) \le 2k \le 2d+d|e|(d+1).
\end{equation}
 It is easy to show that the  positive  integer $d$ satisfying (\ref{eq25}) is unique. Conversely,
for $d$ satisfying (\ref{eq25}), we have unique solution $(m_+, m_-)$ to (\ref{eq13}).

Therefore,
$$c_k(Y, \lambda) =\lim_{\varepsilon \to 0}c_k(Y, \lambda_{\varepsilon}) =\mathcal{A}_{\lambda_{\varepsilon}}(e_+^{m_+} e_-^{m_-})= 2d|e|,   $$
where $d$ is the positive integer  satisfying (\ref{eq25}).

 % Therefore, $i_{L_{\varepsilon}}: ECH_{2k}^{L_{\varepsilon}} (Y, \lambda_{\varepsilon}, 0) \to ECH_{2k}(Y, \lambda_{\varepsilon}, 0)$ is an isomorphism. % Since  $ECH_{2k}^{L} (Y, \lambda_{\varepsilon}, 0) =<e_-^{m_-}e_+^{m_+}>$ for $2d|e| < L\le L_{\varepsilon}$, where $(m_-, m_+)$ is the unique pair such that $gr(e_-^{m_-}e_+^{m_+}) =2k$. We must have $c_{k}(Y, \lambda_{\varepsilon})  = A_{\lambda_{\varepsilon}}(e_-^{m_-}e_+^{m_+}) =2d|e|$.

%Since the ECH capacities is nondecreasing  with respective to $c_{k}(Y, \lambda_{\varepsilon}) \le c_{k_0}(Y, \lambda_{\varepsilon}) <L_{\varepsilon}$ for any $0 \le k \le k_0$.

%By definition, for any $c_k(Y, \lambda_{\varepsilon})  \le 2d|e| L< L_{\varepsilon}$, $i_L: ECH_{2k}^L(Y, \lambda_{\varepsilon}, 0)  \to ECH_{2k}(Y, \lambda_{\varepsilon}, 0)$ is an isomorphism. Therefore, $c_k(Y, \lambda_{\varepsilon}) =2d|e|$.

\end{proof}

\subsection{Torus case}
By definition, $c_k(Y, \lambda) =\infty$ if we cannot find $\sigma \in ECH(Y, \lambda, 0)$ such that $U^k\sigma =[\emptyset]$.  Using computations of P. Ozsv\'{a}th, Z. Szab\'{o} \cite{OS}, and K. Park \cite{Park}, the existence of
such the classes can be guaranteed.
\begin{lemma} \label{lem15}
There exists a sequence of classes $\sigma_{k} \in ECH(Y, \lambda, 0)$ such that $U(\sigma_{k}) =\sigma_{k-1}$ and $U^{k} (\sigma_{k}) =[\emptyset]$.
\end{lemma}
\begin{proof}
By the results of V. Colin, P. Ghiggini and K. Honda \cite{CGH1, CGH2, CGH3}, we know that
$ECH(Y, \Gamma)$ is isomorphic to the Heegaard Floer homology $HF^+(-Y, \mathfrak{s}_{\Gamma})$ as an $\mathbb{F}[U]$-module.

One can see this isomorphism
alternatively by Taubes's isomorphism ``ECH=SWF" \cite{T2, Te2, Te3, Te4, Te5}, and
also by C. Kutluhan, Y-J. Lee, and CH. Taubes’s isomorphism ``SWF=HF" \cite{KLT1, KLT2,KLT3,KLT4,KLT5}. Here ``SWF" stands for the Seiberg-Witten Floer homology and  ``HF" for $HF^+(-Y, \mathfrak{s}_{\Gamma})$.

Thanks to Theorem 5.6 of \cite{OS} and Theorem 4.1.1 of \cite{Park}, we have
 \begin{equation*}
   HF^+(-Y, \mathfrak{s}_0) \cong   \mathbb{F}[[U] \oplus \mathbb{F}[[U] \oplus \mathbb{F}[[U]\oplus \mathbb{F}[[U]
 \end{equation*}
where $ \mathbb{F}[[U] := \mathbb{F}[U^{-1},U]/U\mathbb{F}[U]$. Therefore, $\sigma_k := U^{-k}([\emptyset])$ are well defined. Since $U$ is
degree $-2$ and $\operatorname{gr}([\emptyset]) = 0$, we have $\operatorname{gr}(\sigma_k) = 2k$.

\end{proof}

\begin{proof}[Proof of Theorem \ref{thm4}]
As in the proof of Theorem \ref{thm3}, Lemma \ref{lem15} and Theorem 1.3 of \cite{CHR}
imply  that $c_k(Y, \lambda_{\varepsilon})$ is finite for each $k$. Then for each $k$ and $0 < \varepsilon \ll 1$, we have
$c_k(Y, \lambda_{\varepsilon}) < L_k <L_{\varepsilon}$.

By Lemma \ref{lem14}, we have an ECH generator $\alpha$ such that $c_k(Y, \lambda_{\varepsilon})= \mathcal{A}_{\lambda_{\varepsilon}}(\alpha)$ and $\operatorname{gr}(\alpha) =2k$.
The  ECH generator $\alpha=e_-^{m_-}h_1^{m_1}h_2^{m_2}e_+^{m_+}$ satisfying $\operatorname{gr}(\alpha) =2k$ and $[\alpha]=0$  is
equivalent to
 \begin{equation} \label{eq19}
    \begin{split}
     &\operatorname{gr}(\alpha)=d^2|e| +m_+-m_-=2k\\
     &m_++m_1+m_2+m_-=d|e| \mbox{ and } m_1, m_2 \in \{0, 1\}
    \end{split}
 \end{equation}
  for some $d \in \mathbb{Z}_{\ge1}$.  From the relations (\ref{eq19}), it is easy check that $d$ satisfy
  \begin{equation*}
    d(d-1) |e| \le 2k \mbox{ and }     d(d+1) |e| \ge  2k.
  \end{equation*}
  Solve the inequality $d(d -1) \le  2k/|e| \le d(d + 1)$; we get
  \begin{equation*}
    \sqrt{\frac{2k}{|e|} + \frac{1}{4}} - \frac{1}{2} \le  d \le \sqrt{\frac{2k}{|e|} + \frac{1}{4}} + \frac{1}{2}.
  \end{equation*}
  Let $d_{\max} := \sqrt{\frac{2k}{|e|} + \frac{1}{4}} + \frac{1}{2} $
and $d_{\min}: =    \sqrt{\frac{2k}{|e|} + \frac{1}{4}} - \frac{1}{2}.$   Since
$d$ is an integer, we have $ \lceil  d_{\min}\rceil \le d\le \lfloor d_{\max}\rfloor$, where $\lceil  d_{\min}\rceil$   is the minimal integer that
is greater than or equal to $d_{\min}$,  and $ \lfloor d_{\max}\rfloor$ is the maximal integer  that is no more than $d_{\max}$.

 Assume that $k \ne  \frac{n(n+1)|e|}{2}$ for any $n  \in \mathbb{N}$. Then $d_{\min}$ is not an integer.  Hence, we
have $d_{\min}  =\lfloor  d_{\min}\rfloor +r$  for some $0 < r < 1$.  Note that $d_{\max} =d_{\min} +1$. Then
\begin{equation*}
\lfloor d_{\max}\rfloor =\lfloor \lfloor d_{\min}\rfloor +r + 1 \rfloor =\lfloor d_{\min}\rfloor + 1 = \lceil  d_{\min}\rceil.
\end{equation*}
  As a result, $d= \lfloor d_{\max}\rfloor  = \lceil  d_{\min}\rceil$. Therefore,
  \begin{equation*}
    c_k(Y, \lambda) = \lim_{\varepsilon \to 0} c_k(Y, \lambda_{\varepsilon}) = \lim_{\varepsilon \to 0} 2 \lfloor d_{\max}\rfloor  |e|+O(\varepsilon) = 2 \lfloor d_{\max}\rfloor  |e|.
  \end{equation*}

  Now we consider the case that $k = |e|$. Let $\sigma \in ECH_{2|e|}(Y, \lambda_{\varepsilon}, 0)$ such that $U^{|e|}\sigma= [\emptyset]$. By the index formula (\ref{eq8}), $e_+^{d|e|}$
 and $e_-^{(d+1)|e|}$ are the only ECH generators with
grading $d(d + 1)|e|$. By Lemma \ref{lem12}, the map $U^{|e|}$  is of the form
\begin{equation*}
 U^{|e|}\begin{pmatrix}
e_+^{|e|} \\
e_-^{2|e|}
\end{pmatrix} = \begin{pmatrix}
a & 1 \\
b & c
\end{pmatrix}\begin{pmatrix}
\emptyset\\
e^{|e|}_-
\end{pmatrix},
\end{equation*}
By the computation in Lemma \ref{lem15}, we know that $U^{|e|} : ECH_{2|e|}(Y, \lambda, 0) \to  ECH_0(Y, \lambda, 0)$
is an isomorphism. Thus, we have $b + ac = 1$. If $b = 0$, then $a = c = 1$. If $b = 1$, then
$a = 0$ or $c = 0$. In either cases, we have $U(e_+^{|e|}
+ e_-^{2|e|}) = \emptyset$ or $U(e^{2|e|}_-) = \emptyset$. By (\ref{eq17}),
we have $$c_{|e|}(Y, \lambda) = \lim_{\varepsilon \to 0} c_{|e|}(Y, \lambda_{\varepsilon}) = \lim_{\varepsilon \to 0} \mathcal{A}_{\lambda_{\varepsilon}}(e_-^{2|e|} ) = 4|e|.$$
  \end{proof}

%\begin{lemma}
%Assume that $\Sigma$ is the torus, then one of  the following statements is true:
%\begin{enumerate}
%  \item
%  $\sigma_{2k_0} = i_L(e_-^{M+|e|} +a e_+^M)$;
%  \item
%   $\sigma_{2k_0-2M} = i_L(e_-^{M} +a e_+^{M-|e|})$.
%\end{enumerate}
%\end{lemma}

Shenzhen University
\verb| E-mail address: ghchen@szu.edu.cn|

\end{document}